\documentclass[12pt]{amsart}
\usepackage{amsbsy}
\begin{document}

\newtheorem{theorem}{Theorem}
\newtheorem{proposition}{Proposition}
\newtheorem{lemma}{Lemma}
\newtheorem{corollary}{Corollary}
\newtheorem{definition}{Definition}
\newtheorem{remark}{Remark}
\newcommand{\tex}{\textstyle}
\numberwithin{equation}{section} \numberwithin{theorem}{section}
\numberwithin{proposition}{section} \numberwithin{lemma}{section}
\numberwithin{corollary}{section}
\numberwithin{definition}{section} \numberwithin{remark}{section}
\newcommand{\ren}{\mathbb{R}^N}
\newcommand{\re}{\mathbb{R}}
\newcommand{\n}{\nabla}
\newcommand{\iy}{\infty}
\newcommand{\pa}{\partial}
\newcommand{\fp}{\noindent}
\newcommand{\ms}{\medskip\vskip-.1cm}
\newcommand{\mpb}{\medskip}
\newcommand{\AAA}{{\bf A}}
\newcommand{\BB}{{\bf B}}
\newcommand{\CC}{{\bf C}}
\newcommand{\DD}{{\bf D}}
\newcommand{\EE}{{\bf E}}
\newcommand{\FF}{{\bf F}}
\newcommand{\GG}{{\bf G}}
\newcommand{\oo}{{\mathbf \omega}}
\newcommand{\Am}{{\bf A}_{2m}}
\newcommand{\CCC}{{\mathbf  C}}
\newcommand{\II}{{\mathrm{Im}}\,}
\newcommand{\RR}{{\mathrm{Re}}\,}
\newcommand{\eee}{{\mathrm  e}}
\newcommand{\eb}{{\bf  e}}
\newcommand{\LL}{L^2_\rho(\ren)}
\newcommand{\LLL}{L^2_{\rho^*}(\ren)}
\renewcommand{\a}{\alpha}
\renewcommand{\b}{\beta}
\newcommand{\g}{\gamma}
\newcommand{\G}{\Gamma}
\renewcommand{\d}{\delta}
\newcommand{\D}{\Delta}
\newcommand{\e}{\varepsilon}
\newcommand{\var}{\varphi}
\newcommand{\lll}{\l}
\renewcommand{\l}{\lambda}
\renewcommand{\o}{\omega}
\renewcommand{\O}{\Omega}
\newcommand{\s}{\sigma}
\renewcommand{\t}{\tau}
\renewcommand{\th}{\theta}
\newcommand{\z}{\zeta}
\newcommand{\wx}{\widetilde x}
\newcommand{\wt}{\widetilde t}
\newcommand{\noi}{\noindent}
\newcommand{\uu}{{\bf u}}
\newcommand{\xx}{{\bf x}}
\newcommand{\yy}{{\bf y}}
\newcommand{\zz}{{\bf z}}
\newcommand{\aaa}{{\bf a}}
\newcommand{\cc}{{\bf c}}
\newcommand{\jj}{{\bf j}}
\newcommand{\ggg}{{\bf g}}
\newcommand{\UU}{{\bf U}}
\newcommand{\YY}{{\bf Y}}
\newcommand{\HH}{{\bf H}}
\newcommand{\GGG}{{\bf G}}
\newcommand{\VV}{{\bf V}}
\newcommand{\ww}{{\bf w}}
\newcommand{\vv}{{\bf v}}
\newcommand{\hh}{{\bf h}}
\newcommand{\di}{{\rm div}\,}
\newcommand{\ii}{{\rm i}\,}
\newcommand{\nn}{{\bf  n}}
\newcommand{\inA}{\quad \mbox{in} \quad \ren \times \re_+}
\newcommand{\inB}{\quad \mbox{in} \quad}
\newcommand{\inC}{\quad \mbox{in} \quad \re \times \re_+}
\newcommand{\inD}{\quad \mbox{in} \quad \re}
\newcommand{\forA}{\quad \mbox{for} \quad}
\newcommand{\whereA}{,\quad \mbox{where} \quad}
\newcommand{\asA}{\quad \mbox{as} \quad}
\newcommand{\andA}{\quad \mbox{and} \quad}
\newcommand{\withA}{,\quad \mbox{with} \quad}
\newcommand{\orA}{,\quad \mbox{or} \quad}
\newcommand{\atA}{\quad \mbox{at} \quad}
\newcommand{\onA}{\quad \mbox{on} \quad}
\newcommand{\ef}{\eqref}
\newcommand{\ssk}{\smallskip}
\newcommand{\LongA}{\quad \Longrightarrow \quad}
\def\com#1{\fbox{\parbox{6in}{\texttt{#1}}}}
\def\N{{\mathbb N}}
\def\A{{\cal A}}
\newcommand{\de}{\,d}
\newcommand{\eps}{\varepsilon}
\newcommand{\be}{\begin{equation}}
\newcommand{\ee}{\end{equation}}
\newcommand{\spt}{{\mbox spt}}
\newcommand{\ind}{{\mbox ind}}
\newcommand{\supp}{{\mbox supp}}
\newcommand{\dip}{\displaystyle}
\newcommand{\prt}{\partial}
\renewcommand{\theequation}{\thesection.\arabic{equation}}
\renewcommand{\baselinestretch}{1.1}
\newcommand{\Dm}{(-\D)^m}

\title
{\bf  Boundary characteristic point regularity for Navier--Stokes
equations:  blow-up scaling and Petrovskii-type criterion (a
formal approach)}

\author {V.A.~Galaktionov and V.~Maz'ya}

\address{Department of Mathematical Sciences, University of Bath,
 Bath BA2 7AY, UK}
\email{vag@maths.bath.ac.uk}

\address{Department of Mathematical Sciences, University of Liverpool,
  M$\&$O Building,
Liverpool, L69 3BX, UK \,\, {\rm and} \,\, Department of
Mathematics, Link\"oping University, SE-58183, Link\"oping,
Sweden}
 \email{vlmaz@liv.ac.uk \, and \,  vlmaz@mai.liu.se}

\keywords{Navier--Stokes equations in $\re^3$, backward
paraboloid, characteristic vertex, boundary regularity, blow-up
scaling, boundary layer, solenoidal Hermite polynomials,
eigenfunction expansion, matching,
   Petrovskii's criterion, fourth-order Burnett equations.}

 \subjclass{35K55, 35K40}
\date{\today}


\begin{abstract}

The {\em three-dimensional} (3D) {\em Navier--Stokes equations}
  \be
  \label{0.1}
 \uu_t+ (\uu \cdot \n)\uu=- \n p+\D \uu, \quad {\rm div} \, \uu=0
 \inB Q_0,
  \ee
  where $\uu=[u,v,w]^T$ is the vector field  and $p$ is the
  pressure,
 are considered. Here, $Q_0 \subset \re^3 \times [-1,0)$
   is a smooth domain of a typical {\em backward paraboloid} shape, with the vertex
   $(0,0)$ being its only {\em characteristic point}: the plane $\{t=0\}$ is
   tangent to $\pa Q_0$ at the origin,
 and other characteristics for $t \in [0,-1)$ intersect $\pa Q_0$
   transversely. Dirichlet boundary
   conditions on the lateral boundary $\pa Q_0$ and smooth initial data are
   prescribed:
     \be
     \label{0.2}
     \uu=0 \quad \mbox{on} \quad \pa Q_0, \andA
 \uu(x,-1)=\uu_0(x) \inB Q_0 \cap \{t=-1\} \quad ({\rm div} \, \uu_0=0).
  \ee
Existence, uniqueness, and regularity studies of \ef{0.1} in {\em
non-cylindrical domains} were initiated in the 1960s  in
pioneering works by J.L.~Lions, Sather, Ladyzhenskaya, and
Fujita--Sauer. However, the problem of a {\em characteristic
vertex} regularity remained open.

In this paper, the classic  problem of regularity (in Wiener's
sense) of the vertex  $(0,0)$ for \ef{0.1}, \ef{0.2}
 is considered.
 Petrovskii's famous ``$2\sqrt{\log\, \log}$-criterion" of boundary regularity for the heat equation (1934)
  is shown to apply.
  Namely, after a blow-up scaling
  and a special matching with a boundary
  layer near $\pa Q_0$, the regularity problem reduces
  to
 a 3D
  perturbed nonlinear dynamical system  for the first
 Fourier-type
 coefficients
 of the solutions
  expanded using
 {solenoidal}  Hermite polynomials.
 Finally, this confirms that the nonlinear convection term gets an exponentially decaying factor and is
 then  negligible. Therefore, the regularity of the vertex is
 entirely dependent on the linear terms and hence remains the same
 for  Stokes' and purely parabolic problems.

  Well-posed Burnett equations with the minus bi-Laplacian in \ef{0.1} are also
  discussed.


\end{abstract}

\maketitle

\begin{center}
{\em This work is dedicated  to Professor V.~Lakshmikantham with  great esteem}


\end{center}




\section{Introduction: vertex regularity for the Navier--Stokes equations}
\label{S0}

\subsection{Navier--Stokes equations inside a non-cylindrical backward
paraboloid: first history since 1960s}

 We consider  3D
Navier--Stokes equations (the 3D NSEs)
  \be
  \label{1}
 \uu_t + (\uu \cdot \n)\uu=- \n p + \D \uu, \quad {\rm div} \, \uu=0 \inB Q_0,
   \ee
  where $\uu=[u,v,w]^T(x,t)$ is the  vector field and $p=p(x,t)$ is the corresponding  pressure.

   The NSEs \ef{1} are posed in a smooth {\em non-cylindrical} domain
    $$
    Q_0 \subset \re^3 \times
    [-1,0)
    $$
 of a typical {\em backward paraboloid} shape, with the vertex
   $(0,0)$
     being its only {\em characteristic point}: the plane $\{t=0\}$ is
   tangent to $\pa Q_0$ at the origin. No characteristic points of
   $\pa Q_0$ are assumed to exist for $t \in [-1,0)$, i.e., other characteristics
    $\{t=\t\}$, for any $\t \in [-1,0)$,
    intersect $\pa Q_0$
   transversely, in a natural sense.
 Next,
 the zero Dirichlet boundary
   conditions on the lateral boundary $\pa Q_0$ and smooth initial data at   $t=-1$ are
   prescribed:
     \be
     \label{u01}
     \begin{matrix}
     \uu=0 \quad \mbox{on} \quad \pa Q_0, \andA \ssk\ssk\\
 \uu(x,-1)=\uu_0(x) \inB Q_0 \cap \{t=-1\} \whereA   {\rm div} \, \uu_0=0.
  \end{matrix}
  \ee

The questions of solvability, uniqueness, and regularity for the
Navier--Stokes equations in
  {\em non-cylindrical} (and non-characteristic) domains, i.e., in our case, up to the
  vertex, for $t \le -\d_0<0$,
were actively studied since the 1960s. J.L.~Lions began this study
in 1963;
 see references in his
   classic monograph \cite[Ch.~3]{JLi} concerning elliptic
  regularization--penalization methods; as well as Fujita--Sauer \cite{FujS69, FujS70} (1969)
 as one of the first such study of weak solutions via a penalization.
    Another alternative, as
  was pointed out in \cite[Ch.~3, \S~8.1]{JLi}, is a ``rather careful using" Galerkin methods
with time dependent basis functions; see Sather \cite{Sather63}
(1963).
   In 1968,
Ladyzhenskaya
 \cite{Lad68} proved local existence (global for $N=2$ and for small initial data if $N=3$)
  and uniqueness of strong solutions for time-dependent domains using a
 different method. See Neustupa \cite{Neust09} for more recent
 results, references, and other related problems.

 \ssk

 However, the problem of regularity of a characteristic boundary
point for the NSEs in any dimension $N \ge 2$ was not addressed
elsewhere and remained open. Naturally, in order to proceed with
regularity issues concerning the paraboloid vertex $(0,0)$,
  we have to  assume that a unique smooth bounded
  solution of \ef{1}, \ef{u01} exists in $Q_0$, i.e., with no $L^\iy$-blow-up for
   $t<0$\footnote{But the
  solution is formally allowed to blow-up at the vertex $(0,0^-)$.}.
  In particular, as is well-known (see  \cite{Lad68, Maj02, Temam85}), global
  smooth solutions always exist for sufficiently small initial
  data, so we can directly proceed with  the vertex regularity, at least,
  for this class of solutions.


\subsection{Regularity of the characteristic paraboloid  vertex}
 Thus,
the classic  problem of regularity (in {\em Wiener's sense}, see
\cite{Maz99}) of a boundary characteristic point for the NSEs
problem \ef{1}, \ef{u01}
 is under consideration.

\ssk

\noi{\bf Definition (vertex regularity/irregularity).}
 {\em According to Wiener \cite{Wien24}, the vertex $(0,0)$ of the given backward paraboloid
 $Q_0$ for the NSEs problem \ef{1}, \ef{u01}  is
 {\bf
 regular} if, for any bounded  data $\uu_0(x)$,}
  \be
  \label{2}
  \uu(0,0^-)=0,
   \ee
   {\em and {\bf irregular} otherwise, i.e., at least for some initial data,
   \ef{2} fails.}

\ssk

 The  boundary and other regularity issues for the Navier--Stokes
equations  in $\re^2$ and $\re^3$ have been and remain key and
very popular in modern mathematical literature,
 since J.~Leray's
seminal papers  in 1933-34 \cite{Ler33, Ler34}. Among various
regularity and partial regularity results for the NSEs, the
boundary regularity properties in piecewise smooth or Lipschitz
domains and those with thin channels, or other non-regular domains
(as we will show, such settings are key for our study) always
played a special role. Mentioning Kondratiev's first study of 1967
\cite{Kond67}, we refer to advanced  results, further references,
and  reviews in recent papers \cite{DeurWolf95, Kohr2010, Kweon07,
MazRossNSE09.1, MazRossNSE09.2, MazNSE00, Mit08} and
\cite{Wolf10}. See also \cite{KozMaz01, Mit09} for  the related
linear Stokes problem
\be
 \label{Stokes1}
  \uu_t=- \n p + \D \uu, \quad {\rm div}\, \uu=0 \inB Q_0,
  \quad \uu(0,x)=\uu_0(x) \,\,\,({\rm div}\, \uu_0=0).
   \ee
Concerning compressible flows and other related problems, see a
good survey in \cite{Kweon07}, where 2D NSEs in a polygon domain
with a convex vertex were studied.


Note that, and this is key for us in what follows, J.~Leray in
\cite{Ler33, Ler34} actually posed a deep
 problem on both {\em backward} and {\em forward continuation}
 phenomena, which sound modern and advanced
 nowadays for general nonlinear PDE theory:
 \be
 \label{for1}
 \begin{matrix}
 \mbox{{\bf Leray's blow-up scenario:} \quad
self-similar blow-up as $t \to T^-$ ($t<T$)}\qquad\qquad \ssk\\
 \mbox{and similarity collapse of this
singularity as $t \to T^+$ ($t
> T$)}; \qquad\qquad
 \end{matrix}
 \ee
  see his precise statements and a discussion on these principal issues
  in \cite[\S~2.2]{GalNSE}. In this connection, such ``backward blow-up
  scaling approaches" will be key later  on.

  \ssk


 According to our approach, we deal with a  typical asymptotic problem
of clarifying a generic behaviour of solutions near a ``blow-up
singularity" $(0,0)$ (a ``micro-scale structure" of nonlinear PDEs
involved). Of course, the vertex regularity problem setting
essentially and crucially depends on the {\em a priori} given
shape of the prescribed backward paraboloids, which affects our
methods of matched blow-up expansions. Anyway,  we hope that our
blow-up analysis on shrinking as $t \to 0^-$ subsets
will eventually help to better understand the possible nature of
other plausible  blow-up singularities of the NSEs.  As a common
feature of our blow-up analysis, we will see that a
complete/closed set of vector solenoidal Hermite
polynomials\footnote{These are the same vector polynomials that
occur in the study of multiple zero formation for the Stokes
equations and NSEs, \cite[\S~3]{GalNSE}.} in $\re^3$ can play an
important part.

\ssk

Our main goal here is as follows: using  techniques of blow-up
scaling and matched asymptotic expansions from reaction-diffusion
theory, {\em to show that a Petrovskii's-like criterion of
boundary regularity for the heat equation $(1934)$  does occur for
the NSEs \ef{1}}.



\subsection{Petrovskii's criterion of 1934 for the heat equation: a first discussion}

 For the {\em heat equation},
 \be
 \label{heat1}
 u_t=   \D u \inB Q_0,
  \ee
 the regularity problem of the characteristic vertex $(0,0)$ was {\em optimally} solved for
 dimensions $N=1$ and
 2
  by {\em Ivan Georgievich  Petrovskii} in 1934-35
 \cite{Pet34Cr, Pet35}\footnote{Compare with the dates of Leray's pioneering study, ``1933-34"; almost a perfect  coincidence!},
  who introduced his famous
{\em Petrovskii's regularity criterion}; see \cite{GalPet2m} and
\cite{GalMazf(u)}  for a full history and further developments in
general parabolic theory. This  is the so-called ``$2\sqrt{\log
\log}$--criterion" (see \ef{ph2} and \ef{PCr1} below), which we
are going to
 achieve, at least formally,
   for the 3D NSEs.

\ssk

 The following issues naturally occur:

\ssk

 \noi{\bf (i)} On one hand, Petrovskii's-like criterion can be
 expected, since \ef{1}, similar to \ef{heat1}, is indeed a (vector) parabolic second-order
 equation;

\ssk

 \noi{\bf (ii)} On the other hand, \ef{1} is a {\em  nonlocal} parabolic
 PDE for solenoidal vector fields, and it is not straightforward
 from the beginning
 that this cannot affect the regularity standing; and

 \ssk

 \noi{\bf (iii)} Finally and most essentially, \ef{1} contains
 both linear (the second-order Laplacian) and nonlinear (convective)
  operators, so that the regularity of the vertex $(0,0)$ inevitably will
depend on both, which makes the analysis more difficult.

\ssk

 Note that both the issues {\bf (i)} and {\bf (ii)} apply to
the linear Stokes problem without the quadratic convection
\ef{Stokes1},
 for which our regularity results turn out to be new as well.

 \subsection{Layout of the paper}

In Section \ref{Sect1}, we perform first a {\em blow-up} scaling
near the characteristic vertex $(0,0)$.
 In Sections \ref{S3} and \ref{S3.4}, the main goal is to show how
 the convection term in the NSEs \ef{1}
can affect the regularity conditions by deriving sharp formal
asymptotics of solutions near the characteristic point.  A
necessary and  already well-existing spectral theory involving a
complete set of vector {\em solenoidal Hermite polynomials} as
eigenfunctions of the linear Hermite operator is described in
Appendix A.

For the sake of our regularity study,
 we apply a method of a matched asymptotic (blow-up)
expansion, where a {\em Boundary Layer} behaviour close to the
lateral boundary $\partial Q_0$ (Section \ref{S3}) is matched, as
$t \to 0^-$, with a {\em ``centre subspace behaviour"} in an {\em
Inner Region}, developed in Section \ref{S3.4}. This leads to a
perturbed 3D nonlinear dynamical system for the first Fourier-like
coefficients in the eigenfunction expansions of the vector field
$\uu(x,t)$  via standard solenoidal Hermite polynomials.
 This
approach falls into the scope of typical ideas of asymptotic PDE
theory, which got a full mathematical justification for many
problems of interest. We refer to a most general asymptotic
analysis performed in \cite{KozMaz04}, and also to a number of
complicated blow-up asymptotics in reaction-diffusion theory
\cite{AMGV}. According to the classification in \cite{KozMaz04},
our matched blow-up approach corresponds to perturbed
three-dimensional dynamical systems, i.e., to a rather
not-that-advanced case being however a constructive one that has
given a number of new asymptotic/regularity results.
 We propose a final, more general discussion of various (boundary and interior) blow-up singularities for the NSEs in Section \ref{SFi}.


\ssk

In addition, for more clear expressing  our regularity techniques
 and their applicability in general  PDE theory, we develop
in  Appendix B (the C one contains the corresponding Hermite
spectral analysis) at the paper end,
 as a natural extension,
 a similar regularity analysis of the {\em well-posed Burnett
equations} in $Q_0$ with zero Dirichlet conditions
 \be
 \label{1B}
  \left\{
  \begin{matrix}
 \uu_t+ (\uu \cdot \n)\uu=-\n p  - \D^2 \uu, \quad{\rm div} \, \uu=0 \inB Q_0, \qquad\qquad\qquad\ssk\ssk\\
 \uu=\n\uu \cdot \nn=0 \quad \mbox{on} \quad \pa Q_0,
 \quad \uu(x,-1)=\uu_0(x), \qquad\qquad\qquad\quad\,\,\,
  \end{matrix}
  \right.
  \ee
  where $\nn$ denotes the unit inward normal vector to $\pa Q_0
  \cap \{t\}$.
Here we have the {\em bi-harmonic} diffusion operator $-\D^2 \uu$  on the
right-hand side of the $\uu$-equation.
 It turns out that our general scheme of the boundary regularity
 analysis can be applied; however, with harder asymptotics  and  more formal nature of
 the final difficult estimates.

\ssk

Concerning other problems and techniques of modern regularity
theory, we  refer to monographs \cite{Gris85, KMR1, KozMaz01,
Maz2010, MazSol2010} and \cite{KozMaz05},
\cite{MayMaz09}--\cite{MS04}  as an update guide to elliptic
regularity theory including higher-order equations, as well as to
references/results in
 \cite{Gris92, Kweon07, Land69, Lanc73, Evans82, Wat97} and \cite{GalPet2m, GalMazf(u)}
  for linear and semilinear parabolic PDEs.


 \section{First blow-up scaling:  Sturm's backward scaling
 variable, paraboloid geometry, and   the Cauchy problem}
  \label{Sect1}


\subsection{First blow-up scaling: exponentially small perturbations of a rescaled parabolic flow}

  We perform
   {\em blow-up scaling} in (\ref{1})
    for the regularity analysis
   \be
   \label{ll2}
   \tex{
   \uu(x,t)=  \vv (y,\t), \quad y = \frac x{\sqrt{-t}}, \quad \t=
   -\ln(-t): (-1,0) \to \re_+,
 }
 \ee
 where $y$ is, indeed, {\em Sturm's backward rescaled  variable} introduced him in
 1836 \cite{St36} in the study of zero sets of solutions of linear parabolic
 equations such as \ef{heat1} for $N=1$.

Thus,
 scaling \ef{ll2} yields
 the following {\bf exponentially perturbed}  rescaled
 equation\footnote{Here, \ef{S1}
 {\em is not} the rescaled one \ef{Full1} written in Leray's
 variable \ef{ll2Ler}, a difference to be discussed.}:
 \be
  \label{S1}
   \begin{matrix}
    \vv _\t + \eee^{-\frac \t 2}\, (\vv
\cdot \n)\vv= - \eee^{-\frac \t 2}\, \n p+
     \BB^* \vv, \quad {\rm div} \, \vv=0    \inB \hat Q_0,
     \ssk\ssk\\
     \mbox{where} \quad  \vv=[v^1, v^2, v^3]^T \andA
  \BB^*=\D   -
 \frac 12 \, y \cdot \n
  \end{matrix}
  \ee
  is Hermite's classic symmetric (self-adjoint) operator, \cite[p.~48]{BS}.


 \subsection{Backward paraboloid geometry and a slow growing factor $\var(\t)$}

According to Petrovskii \cite{Pet34Cr, Pet35}, the backward
paraboloid $\pa Q_0$ will be defined as follows: it is a
perturbation of the standard {\em fundamental} backward one,
 \be
 \label{ph1}
  \tex{
 S(t)= \pa Q_0\cap\{t\}: \quad q_0(x) \equiv \sqrt{ \sum_{i=1}^3 a_i x_i^2}=
  (-t)^{\frac 12}\, \var(\t),
 }
  \ee
  where $a_i \in (0,1]$ are normalized constants.
We can treat more general {\em convex} paraboloids, but, for
simplicity, will restrict to the basic ones as in \ef{ph1}, which
are also of a challenge. For difficult estimates to follow, we
will use {\em radially symmetric} paraboloids with $a_i=1$:
 \be
 \label{rad11}
 \tex{
 S(t)= \pa Q_0\cap\{t\}: \quad q_0(x) \equiv |x|
 = (-t)^{\frac 12}\, \var(\t).
 }
  \ee

 In \ef{ph1} and \ef{rad11},  $\var(\t)>0$ is a slow growing function
 satisfying
 \be
 \label{vv1}
  \tex{
 \var(\t) \to +\iy, \quad \var'(\t)>0,  \quad \var'(\t) \to 0, \andA \frac {\var'(\t)}{\var(\t)} \to 0 \asA \t \to
 +\iy.
  }
  \ee
 Moreover, as a sharper characterization of the above class of
 {\em slow growing functions}, we use the following criterion:
  \be
  \label{vv2}
   \tex{
   \big( \frac{\var(\t)}{\var'(\t)}\big)' \to \iy \asA \t \to
   +\iy.
    }
    \ee
    This is a typical condition in blow-up analysis distinguishing classes of
    exponential (the limit is 0), power-like (a constant $\not =
    0$), and slow-growing functions. See \cite[pp.~390-400]{SGKM},
    where in Lemma 1 on p.~400, extra properties of slow-growing
    functions \ef{vv2} are proved. For instance, one can derive
    the following comparison of such $\var(\t)$ with any power:
     \be
     \label{al1}
     \mbox{for any $\a>0$}, \quad \var(\t) \ll \t^\a \andA
     \var'(\t) \ll \t^{\a-1} \forA \t \gg 1.
     \ee
 Such estimates are useful
 in  evaluating
  perturbation terms in the rescaled equations.

\ssk

In Petrovskii's criterion for the heat equation \ef{heat1}, for  any $N \ge 1$,
 the ``almost
optimal"
 function, satisfying  \ef{vv1}, \ef{vv2} and delivering a regular vertex $(0,0)$, is
  \be
  \label{ph2}
   \fbox{$
   \var_*(\t) = 2\sqrt {\ln \t} \asA \t \to + \iy.
    $}
   \ee
 Replacing this {\em fundamental constant} ``2" by $2+\e$, for an arbitrarily
 small constant $\e >0$, makes $(0,0)$ irregular for the heat equation \ef{heat1}.
In the  general case of arbitrary $\var(\t)$, Petrovskii's
criterion, for our  $N=3$, for the radially symmetric paraboloid  \ef{rad11}
 reads\footnote{For $Q_0 \subset \ren \times [-1,0)$, the
multiplier $\var^3$ in this integral criterion is replaced by
$\var^N$.}
 \be
 \label{PCr1}
 \fbox{$
 \tex{
  \int\limits^\iy \var^3(\t)\, {\mathrm e}^{-\var^2(\t)/4} \,
  {\mathrm d} \t = \iy.
  }
  $}
  \ee
 These  dependencies  will be compared with those obtained for the
 NSEs \ef{1}.

\ssk

Thus, the monotone positive function $\var(\t)$ in \ef{ph1}
  is assumed to determine a sharp behaviour of the
boundary of $Q_0$ near the  vertex $(0,0^-)$ to guarantee
its regularity. It follows that the rescaled equation \ef{S1} is
set in an expanding rescaled domain
  \be
 \label{ph1N}
  \tex{
 \hat S(\t) \equiv  \pa \hat Q_0\cap\{\t\}: \quad  q_0(y) \equiv \sqrt{ \sum_{i=1}^3 a_i y_i^2}=
   \var(\t) \to +\iy
 \asA \t \to +\iy.
 }
  \ee
  By $\nn$ we denote the inward unit normal to $\hat S(\t)$.
In the limit $\t= +\iy$, we arrive at the equation \ef{S1} in the
whole space $\re^3$, which requires some spectral theory (Appendix
A).

\subsection{Towards the Cauchy problem}

In {\em Inner Region} (see Section \ref{S3} for details),
described by compact subsets in the variable $y$ in \ef{ll2}, we
deal with the original rescaled problem
 \ef{S1} in the unboundently expanding domains \ef{ph1N}.
As usual and customary in potential and general PDE theory (see
e.g., Vladimirov  \cite{Vlad72}), it is  convenient to consider
the NSEs in whole space $\re^3 \times [-1,0)$. Note that, in the
study of the NSEs in  non-cylindrical domains, Fujita--Sauer
\cite{FujS69, FujS70} also extended the problem to $\re^3$ by
introducing a {\em  strong absorption term} $-n \, \uu$ in the
complementary domain on the right-hand side in \ef{1} and passing
to the limit $n \to +\iy$. Then, in view of the control of the
total kinetic $L^2$-energy, this ``regularized" solution
$\uu=\uu^{(n)} \to 0$ as $n \to +\iy$ outside the given
non-cylindrical domain, and hence, in the limit, the zero
Dirichlet conditions on the boundary are restored.

 Thus,  we extend
 $\vv(y,\t)$ by 0 beyond the boundary, i.e., set
  \be
  \label{a1}
 \hat \vv(y,\t) = \vv(y,\t)H(\var(\t)-q_0(y))=
  \left\{
   \begin{matrix}
   \vv(y,\t) \forA 0 \le q_0(y) < \var(\t), \\
\,\, 0  \,\,\forA \,\,  q_0(y) \ge  \var(\t),
\end{matrix}
 \right.
 \ee
 where $H(\cdot)$ is the Heaviside function.
Then, since $\vv=0$ on the lateral boundary $\hat S(\t)=\{q_0(y)=
\var(\t)\}$, one can check that, in the sense of distributions
(see e.g., \cite[\S~6.5]{Vlad72}),
 \be
 \label{a2}
  \tex{
 \hat \vv_\t= \vv_\t H, \quad \n \hat \vv= \n \vv H, \quad
 \D \hat  \vv= \D \vv H + \frac {\pa  \vv}{\pa \nn}\,
 \d_{\hat S(\t)},
 }
 \ee
 where   a single-layer potential with the density $\mu= \frac {\pa  \vv}{\pa \nn}$
  acts as follows: for any $\phi \in C_0^\iy$,
  \be
  \label{mat1}
  \tex{
  \langle  \mu \d_{\hat S(\t)}, \phi \rangle= \int_{\hat S(\t)}
  \mu \phi\, {\mathrm d}s.
  }
  \ee

  In order to avoid a pressure trace on $\hat S(\t)$, we perform a
  {\em continuous}  extension of $p(y,\t)$ by solving the
  Laplace equation in the outer domain:
   \be
   \label{ppp1}
    \D \hat p=0 \inB \re^3\setminus (\overline{\hat Q_0(\t)\cap\{\t\}}),
    \quad \hat p=p \,\,\, \mbox{on} \,\,\, \hat S(\t).
     \ee
This outer Dirichlet problem is known to admit a unique solution
$\hat p(y,\t)$ vanishing at infinity, \cite[\S~28]{Vlad72}, so
that we are given the continuous pressure $\hat p(y,\t)$ defined
in the whole $\re^3 \times \re_+$. It then follows that
 \be
 \label{ppp2}
  \begin{matrix}
  \n \hat p= \n p \,\,\,\mbox{in} \,\,\, \overline{\hat Q_0(\t)\cap\{\t\}} \quad (\mbox{since the jump is zero:}
  \,\,\,[\hat p]_{\hat S(\t)}=0), \,\,\, \mbox{and} \qquad \ssk\ssk\\
\n \hat p \,\,\, \mbox{is div-free in $\re^3$ as a
 distribution}: \,\,\, \langle \n \hat p,\phi \rangle=0, \,\,
 \forall\phi \in C_0^\infty, \,\, {\rm div}\, \phi=0. \qquad
 \end{matrix}
  \ee

Thus, $\{\hat \vv, \hat p\}$ satisfies the following Cauchy
problem in $\re^3 \times \re_+$:
 \be
 \label{a3}
  \tex{
  \hat \vv_\t + \eee^{-\frac \t 2}\, (\hat \vv\cdot \n) \hat \vv =
  - \eee^{-\frac \t 2}\,\n \hat p+ \BB^* \hat \vv
   -  \frac {\pa  \vv}{\pa \nn}\,
 \d_{\hat S(\t)}, \quad {\rm div} \, \hat \vv=0.
}
  \ee
 Hence, we obtain a single perturbation term on the right-hand side
expressed in terms of a simple layer potential with the prescribed
density on the surface \ef{ph1N}, changing with the time $\t$.
 Clearly,  various  linear and nonlinear ``interactions" of all these
  and other operators in \ef{a3}  will
 define regularity of the
 vertex.

\subsection{The Cauchy problem in Leray's  nonlocal setting}

Using Leray's nonlocal formulation of NSEs \cite[p.~32]{Maj02}, we
next apply to the Cauchy problem \ef{a3} the operator
 $$
 {\mathbb
P}= I - \n \D^{-1}(\n \cdot
   I) \quad (\| \mathbb{P} \|=1)
    $$
  being the
  Leray--Hopf projector of $(L^2(\re^3))^3$ onto  the subspace
  $\{{\bf w} \in (L^2)^3: \,\, {\rm div}\, {\bf w}=0\}$ of solenoidal vector
   fields.
 Let us note another representation of the projector ${\mathbb P}$
 therein:
   $$
   \mathbb{P} \ww=[v_1-R_1 \s,v_2-R_2 \s,v_3-R_3 \s]^T \whereA \s=R_1 w_1+R_2 w_2+R_3
   w_3,
   $$
   and $R_j$ are the Riesz transforms, with symbols
   $\frac {\xi_j}{|\xi|}$.
 Using the fundamental
 solution of $\D$ in $\ren$, $N \ge 3$, and denoting by
 $\s_N$  the surface area of the unit ball $B_1 \subset \ren$,
  \be
  \label{FF55}
   \tex{
   b_N(y)= - \frac 1{(N-2)\s_N} \, \frac 1{|y|^{N-2}} \whereA \s_N= \frac { 2
   \pi^{ N/2}}{\Gamma( N/2)} \quad(\s_3=4 \pi),
    }
    \ee
the projection of the convective term reads:
 \be
 \label{HHH21}
  \begin{matrix}
- {\mathbb P}(\vv \cdot \n)\vv = - (\vv  \cdot \n )\vv  + C_3
\int\limits_{\re^3} \frac {y-z}{|y-z|^3}\,\, {\rm tr} (\n \vv
(z,\t))^2\, {\mathrm d}z,\qquad \ssk \\
 \mbox{where} \quad
 {\rm tr} (\n
\vv (z,\t))^2= \sum_{(i,j)} \,   v_{z_j}^i  v_{z_i}^j \andA C_N=
\frac 1{\s_N}>0 \,\,\,\big(C_3= \frac 1{4 \pi} \big).\quad\qquad
 \end{matrix}
 \ee
 This is a more convenient form for some estimates.



\ssk

Using the projection ${\mathbb P}$ eliminates the pressure term
$\n \hat p$ in \ef{a3}, and  we obtain the Cauchy problem for the
following perturbed nonlocal parabolic equation for $\hat \vv$:
\be
 \label{a3pp}
 \fbox{$
   \tex{
  \hat \vv_\t   =
   \BB^* \hat \vv - \eee^{-\frac \t 2}\,{\mathbb P}\, (\hat \vv\cdot \n) \hat \vv
   -
     \frac {\pa  \vv}{\pa \nn}\,
 \d_{\hat S(\t)}
  \inB \re^3 \times \re_+.
}
 $}
  \ee
 Since, by construction, the last term {\em is solenoidal}, we
 have omitted the projection ${\mathbb P}$ therein.
 We recall again that
  local existence and uniqueness of a  classic solution
 $\hat \vv(y,\t)$ of \ef{a3pp} are guaranteed by known local regularity properties
 of the NSEs. Moreover, for any sufficiently small data
 $\vv_0$, solutions of \ef{a3pp} are well defined for all $\t \in
 \re_+$, i.e., up to the boundary blow-up moment $t=0^-$
 ($\t=+\iy$). For other solutions, in general, we assume that
 $\hat \vv(y,\t)$ are well defined and do not blow-up at a finite
 $\t>0$, so we need to study their behaviour as $\t \to +\iy$.

 It then follows that Wiener's regularity of the vertex $(0,0)$ is equivalent to the
 following:
  \be
  \label{Poi6}
   \fbox{$
 \begin{matrix}
 \mbox{$0$ is {\em globally asymptotically stable} for $(\ref{a3pp})$,
 i.e.,}\qquad\qquad\qquad\qquad\qquad\qquad\ssk\\
  \mbox{
  any such solution of (\ref{a3pp}) satisfies \,$\hat \vv(y,\t) \to
  0$ as $\t \to +\iy$ uniformly in $y$.}
  \end{matrix}
  $}
  \ee



\subsection{A full pressure representation
 in $\re^3 \times \re_+$}

As usual \cite[p.~30]{Maj02}, once the vector field $\hat \vv$ has
been obtained from \ef{a3pp},
the pressure  is then given by  the corresponding Poisson
equation.
 Since it is slightly  more technical to get it from the rescaled equation
 \ef{a3pp} containing extra operators, we
consider first  the original (non-rescaled) problem for $\hat \uu=
\uu H(\cdot)$, where the Heaviside function $H$ is concentrated on
$Q_0(t) \cap \{t\}$, and construct a harmonic extension $ \hat p$
as in \ef{ppp1} for $S(t)$. This equation, which will be also in
use, is
  \be
  \label{uu1}
   \tex{
  \hat \uu_t +(\hat \uu \cdot \n) \hat \uu= - \n \hat p +
  \D \hat \uu - \frac {\pa \uu}{\pa \nn}\, \d_{S(t)},
  }
   \ee
where we keep the same notation and arguments as in \ef{a3pp}.
Then, taking ${\rm div}$ \cite[p.~30]{Maj02}, we obtain the following equation with
two extra densities of special potentials:
 \be
 \label{uu2}
  \tex{
 - \D_x \hat  p= {\rm tr}\,(\n_x \hat \uu)^2 +
 {\rm div}_x\,(
 \frac {\pa \uu}{\pa \nn_x}\, \d_{S(t)})  - [\n_x p]_{S(t)}\, \d_{S(t)}.
 }
  \ee
Observe that the jump of the pressure gradient $[\n_x p]_{S(t)}$
enters the second density, that makes this elliptic problem more
nonlocal and hence more difficult.

After scaling in  \ef{ll2}, i.e., setting $x=y \sqrt{-t}$ when
approaching the vertex $(0,0)$, we have from \ef{uu2}, taking into
account that $\d_{S(t)} = (-t) \, \d_{\hat S(\t)}$,
 \be
 \label{Poi1}
  \tex{
  -\D_y \hat p = {\rm tr}\, ( \n_y \hat \vv)^2
  +
\eee^{-\t} \, {\rm div}_y\,(
 \frac {\pa \vv}{\pa \nn_y}\, \d_{\hat S(\t)})
  -\eee^{-\frac {3\t} 2}\,  [\n_y p]_{\hat S(\t)}\, \d_{\hat S(\t)}
    \inB \re^3.
   }
   \ee
 Therefore,  this  {\em
nonlocal} problem for $\hat p$  is reduced to a Fredholm linear
integral equation (with a positive kernel) of the second kind,
 \be
 \label{Poi2}
  \tex{
\hat p= (-\D_y)^{-1}\big[ {\rm tr}\, ( \n_y \hat \vv)^2
 +
\eee^{-\t} \, {\rm div}_y\,(
 \frac {\pa \vv}{\pa \nn_y}\, \d_{\hat S(\t)})
  -\eee^{-\frac {3\t} 2}\,  [\n_y p]_{\hat S(\t)}\, \d_{\hat S(\t)}
   \big],
   }
   \ee
   where $(-\D_y)^{-1}$ is defined by the convolution with the fundamental solution
   \ef{FF55}, $N=3$.  It follows that the behaviour of $p$ as $\t \to +\iy$ is
 dependent on the trace of its gradient on the expanded boundary $\hat S(\t)$.
Fortunately, as $\t \to +\iy$, the jump of the gradient $[\n_y
p]_{\hat S(\t)}$  on the right-hand side of \ef{Poi2} has an
exponentially small influence on $\hat p$, so that
\ef{Poi2}/\ef{Poi1} are  ``almost" standard integral/elliptic
equations with good positive kernels/Laplacian. However, we will
need  to justify that, nevertheless, the corresponding densities
of these single and ``double-layer-type" potentials do not get
exponentially large, thus undermining their exponentially small
factors in front of them. This can be done {\em a posteriori},
when the independent rescaled $\hat \vv$-problem \ef{a3pp} has been
solved.

\ssk

   On the other hand, introducing in \ef{Poi1} the following
   integral operator\footnote{For $\t \gg 1$, this is just an asymptotically small
    perturbation of the Laplacian; though proving that $p$ on $\hat S(\t)$ does not grow exponentially fast
    is a part of the problem.}:
    \be
    \label{Poi3}
    {\mathbb M}(\t)p \equiv
  -\D_y \hat p +\eee^{-\frac {3\t} 2}\, [\n_y p]_{\hat S(\t)} \d_{\hat
  S(\t)},
  \ee
  the pressure is given by
   \be
   \label{Poi4}
   \tex{
   \hat p= {\mathbb M}^{-1}(\t) \big(
   {\rm tr}\, ( \n_y \hat \vv)^2
 + \eee^{-\t} \, {\rm div}_y\,(
 \frac {\pa \vv}{\pa \nn_y}\, \d_{\hat S(\t)})\big).
   }
   \ee
Indeed, this is the equivalent pressure representation, since
${\mathbb M}^{-1}(\t)$ is well defined by the local well-posedness
of the NSEs in bounded domains with smooth non-characteristic
boundaries.


 \section{Boundary layer expansion close to $\pa Q_0$}
 \label{S3}



\subsection{Two region expansion}

As we have mentioned,
 by the divergence of $\var(\t)$ in \ef{vv1},
  sharp asymptotics of solutions close to the vertex $(0,0)$
  will
 essentially depend on the spectral properties of the linear
 operator $\BB^*$ in the whole  space $\re^3$ (see Appendix A), as well as on the
 nonlinear projected convective term.
 This ``interaction" between  linear
 and nonlinear operators in the NSEs, together with the paraboloid shape,  are key for us.

 Studying asymptotics of solutions of the rescaled problem \ef{S1},
 as rather often occurs in  difficult blow-up expansions in nonlinear PDE theory, this singularity
 problem is solved by matching of expansions of solutions in two regions:

 \ssk

 (i) In an {\em Inner Region}, which is situated around the origin $y=0$,  and

 \ssk

 (ii) In a {\em Boundary Region} close to the boundary \ef{ph1N}, where a {\em  Boundary
 Layer} occurs.

 \ssk

 \noi In other words, we show that generic behaviour of solutions
 of the NSEs in shrinking neighbourhoods of the paraboloid vertex $(0,0)$
 {\em is not of any self-similar form}, and hence gets more
 complicated and demands novel non-group-similarity techniques to
 detect.

  Actually, such a two-region structure (i)--(ii) above, with the
 asymptotics specified below, defines the class of generic solutions under
 consideration.
We begin with a simpler analysis in the Boundary Region (ii).

 \subsection{Boundary layer (BL) variables and perturbed equation}
   \label{S3.3}

 Sufficiently close to the lateral boundary of $Q_0$, it is
 natural to introduce the next rescaled variables
  \be
  \label{z1}
   \tex{
   z= \frac y{\var(\t)} \andA \hat \vv(y,\t)=\ww(z,\t).
   }
   \ee
This makes the corresponding rescaled paraboloid \ef{ph1N} fixed:
 \be
 \label{QQ1}
  \tex{
   \tilde S:
  \quad \sqrt{\sum_{i=1}^3 a_i z_i^2}=1.
   }
   \ee
The rescaled vector field $\ww$ now solves a perturbed equation:
 \be
 \label{z1Eq}
  \tex{
 \ww_\t=
 \frac 1{\var^2} \, \D_z \ww - \frac 12\, z \cdot \n_z \ww +
    \frac {\var'}\var \, z \cdot \n_z \ww - \frac 1 \var \, \eee^{-\frac \t 2}\,
     {\mathbb P}(\ww \cdot \n_z)\ww
   - 
    \var \,
      \frac {\pa  \ww}{\pa \nn}\,
 \d_{\tilde S}.
 }
    \ee

Let us introduce the BL-variables: fixing a point $z_0 \in \tilde
S$ on the boundary \ef{QQ1},  we set
  \be
  \label{z2}
  \tex{
   \xi= \var^{2}(\t)(z_0-z) \equiv \var(\var z_0-y), \quad \var^{2}(\t) {\mathrm d} \t={\mathrm d}s,
    \quad \ww(z,\t)= \rho(s) \ggg(\xi,s),
  }
  \ee
  where $\rho(s)=[\rho^1(s), \rho^2(s), \rho^3(s)]^T \in \re^3$ for $s \gg 1$ is an unknown
    scaling slow varying/decaying
   (in the same natural sense,
 associated with \ef{vv2})
   time-factor depending on the
  function $\var(\t)$. Cf. e.g.,   $\sim \frac 1{\var(\t)}$ as a clue.
As usual, this $\rho$-scaling is chosen to get uniformly bounded
rescaled solutions, so, we naturally assume that, for each
component,
 \be
 \label{g12}
  \tex{
  \sup_{\xi} \,g^j(\xi,s)=1 \quad \mbox{for all} \quad s \gg 1,
  \quad j=1,\, 2,\,3.
  }
  \ee

On substitution into the PDE in \ef{z1Eq}, we obtain the following
 perturbation of a linear uniformly parabolic equation:
 \be
 \label{z3}
  \begin{matrix}
  \ggg_s= \AAA \ggg - \frac 12\,  \frac 1{\var^{2}} \, \xi \cdot \n_\xi \ggg - \frac
  {\var'_\t}{\var}\, \big(z_0- \frac \xi {\var^{2}}\big)\cdot \n_\xi \ggg
-  \frac{2\var'_\t}{ \var^{3}} \,
  \xi \cdot \n_\xi \ggg -  \frac {\rho'_s}{\rho} \,
  \ggg
   \ssk\ssk\\
 - \frac \rho{\var} \, \eee^{-\frac \t 2} \, {\mathbb P} (\ggg \cdot \n_\xi)\ggg
 - \frac 1{\var^3}\,
  \frac {\pa  \ggg}{\pa \nn_\xi}\,
 \d_{S_\xi(\t)},\,\,\,
   \mbox{where} \,\,\,  \AAA \ggg=  \D_\xi \ggg + \frac 12\, z_0 \cdot \n_\xi
   \ggg,
    \end{matrix}
   \ee
and $S_\xi(\t)$ is the boundary \ef{QQ1} expressed in terms of the
BL-variable $\xi$ in \ef{z2}, so that
 $$
  \tex{
 \d_{S_\xi(\t)}= \frac 1{\var^3(\t)}\, \d_{\tilde S}.
 }
 $$

   As usual in boundary layer theory, the BL-scaling \ef{z2} means that
we  are looking for a generic
 pattern of the behaviour described by the perturbed equation \ef{z3} on  compact
 subsets, shrinking (focusing) to a fixed $z_0$
 on the lateral boundary,
  \be
  \label{z4}
   \tex{
  |\xi| = o\big(\frac 1{\var^{2}(\t)}\big) \to 0
   \LongA |z-z_0| = o\big(\frac 1{\var^{4}(\t)}\big) \to 0 \asA \t \to +
   \iy.
   }
   \ee

Thus, in \ef{z3}, we arrive at a linear uniformly parabolic
equation perturbed by a number of linear and nonlinear terms,
which, under given and other special hypothesis to be specified,
are {\em asymptotically small}. Indeed, on the space-time compact
subsets \ef{z4}, the second term on the right-hand side of \ef{z3}
becomes asymptotically small, while all the other linear ones even
smaller in view of the slow growth/decay assumptions such as
\ef{vv2} for $\var(\t)$ and $\rho(s)$.
 In
particular,
 the rescaled nonlinear convective term in \ef{z3}
 is  asymptotically small on bounded rescaled vector fields $\ggg$
  in view of an exponentially decaying
 factor and  by the hypotheses
 \ef{vv1}.

However, the last  term in \ef{z3}, given by a density of a simple
layer potential, requires a special treatment. Indeed, this
density depends upon  a still unknown gradient of $\ww= \rho \,
\ggg$ on the boundary, which is under scrutiny in the present
BL-analysis. However, the nature of our BL-scaling assumes that we
deal with uniformly bounded rescaled function $\ggg$, on which the
last term is asymptotically small and is of order $ \sim
\frac 1{\var(\t)} \to 0$ as $\t \to
+\iy$.





\ssk

The BL-representation \ef{z2}, by using the rescaling and
\ef{g12}, naturally leads to
the following asymptotic behaviour at infinity:
 \be
 \label{z5}
  \tex{
   \lim_{s \to +\iy} g^j(\xi,s) \to 1 \asA \xi \to  \iy,
   }
  \ee
  where all the derivatives also vanish by the standard interior parabolic regularity.
Actually, the nature of the BL-scaling \ef{z2} near the point $z_0
\in \tilde S$ implies that, asymptotically, the limit problem
becomes one-dimensional, depending on the space variable
 \be
 \label{eta1}
 \eta= \xi \cdot \nn,
  \ee
  where $\nn$ is the unit inward normal to the smooth boundary $\tilde S$
 in \ef{QQ1}. Therefore, the limit $\xi \to \iy$ in \ef{z5} should
 be also understood in the sense of the single variable \ef{eta1}.
 This essentially simplifies the BL-structure to appear.

 Moreover, since according to the BL-scaling \ef{z4}, as $s \to +\iy$, the rescaled
 solution becomes constant (see \ef{z5}) a.e. and hence
 solenoidal, {\em we do not need to require the limit BL-profile to be
 solenoidal as well}.
 Moreover, we will see that, in the actual boundary layer,
 the BL-asymptotic is ``almost" solenoidal, up to an exponentially small perturbation.



 \subsection{Passing to the limit: generic solutions}

   Thus, we arrive at
  the problem of passing to the limit as $s \to + \iy$ in the
  problem \ef{z3}, \ef{z5}. Since, by the definition in \ef{z2},
  the rescaled orbit $\{\ggg(s),\,\, s>0\}$ is uniformly bounded, by
  classic parabolic interior theory \cite{Fr, EidSys, Eid}, one can pass to the
  limit in \ef{z3} along a subsequence $\{s_k\} \to +\iy$. Namely,
   we have that, uniformly on compact subsets defined in
  \ef{z4}, as $k \to \iy$,
   \be
   \label{z6}
   \ggg(s_k+s) \to \hh(s) \whereA \hh_s=\AAA \hh, \quad \hh=0
   \,\,\,\mbox{at} \,\,\,\eta=0, \quad h^j|_{\eta=+\iy}=1.
    \ee

Consider this one-dimensional {\em limit} (at $s=+\iy$) {\em
equation} obtained from \ef{z3}:
 \be
  \label{z7}
   \tex{
  \hh_s= \AAA \hh \equiv  \hh_{\eta\eta}+ \frac 12\, \hh_\eta \inB \re_+ \times \re_+,
  \quad \hh(0,s)=0, \quad
  h^j(+\iy,s)=1.
  }
   \ee
It is a  linear parabolic PDE in the unbounded domain $\re_+$,
governed by the operator
 $\AAA$ admitting a standard symmetric representation in a weighted space.
 Namely, we have:

 \begin{proposition}
  \label{Pr.91}
  {\rm (i)} \ef{z7} is a gradient system in a weighted
  $(L^2)^3$-space, and

  {\rm (ii)} for bounded orbits,  the $\o$-limit set $\O_0$ of \ef{z7} consists of a unique
  stationary profile
  \be
 \label{z12}
  \tex{
  g_0^j(\xi)=1-{\mathrm e}^{- \eta/2}, \quad j=1,2,3,
  }
 \ee
 and $\O_0$ is uniformly stable in the Lyapunov sense in a
 weighted $(L^2)^3$-space.
\end{proposition}

\noi{\em Proof.} (i) As a second-order  equation, \ef{z7}
 can be written in the symmetric form
  \be
  \label{z71}
  {\mathrm e}^{\eta/2} \hh_s= ({\mathrm e}^{\eta/2} \hh_\eta)_\eta,
   \ee
   and hence admits multiplication by $\hh_s$ in $(L^2)^3$ that yields
   a monotone Lyapunov function:
    \be
    \label{z72}
     \tex{
    \frac 12 \, \frac {\mathrm d}{{\mathrm d}s} \int {\mathrm
    e}^{\eta/2}(h^j_\eta)^2=- \int {\mathrm e}^{\eta/2}(h^j_s)^2 \le 0.
    }
    \ee
    Note that, regardless that the weight $\eee^{\eta/2}$ is
    exponentially growing as $\eta \to +\iy$, on the limit profile
    \ef{z12}, all the functionals in \ef{z72} are well defined.
In other words, the problem \ef{z3} is a perturbed {\em gradient
system}, that makes much easier to pass to the limit $s \to +\iy$
by using power tools of gradient system theory; see e.g., Hale
\cite{Hale}.

(ii) For a given bounded orbit $\{h(s)\}$, denote
$h^j(s)=g^j_0+w^j(s)$, so that $\ww(s)$ solves the same equation
\ef{z71}.
 Multiplying by $\ww(s)$ in $(L^2)^3$ yields
  \be
  \label{z73}
   \tex{
    \frac 12 \, \frac{{\mathrm d}}{{\mathrm d}s} \, \int{\mathrm
    e}^{\eta /2} (w^j)^2 \, {\mathrm d} \eta=- \int{\mathrm
    e}^{ \eta /2}(w^j_\eta)^2 \, {\mathrm d} \eta <0
    }
    \ee
for any nontrivial solution, whence the  uniform stability
(contractivity) property. \quad $\qed$

 \ssk

Finally, we state the main stabilization result in the boundary
layer, which also establishes the actual class of {\em generic}
solutions we are dealing with.

\begin{proposition}
 \label{Pr.Gen1}
   Under specified above assumptions and hypotheses, there exists a class of solutions of
   the perturbed equation \ef{z3}, for
 which, in a weighted $(L^2)^3$-space and uniformly on compact
 subsets,
  \be
  \label{z10}
  g^j(\xi,s) \to g^j_0(\xi) \asA s \to +\iy \quad (j=1,2,3).
  \ee

 \end{proposition}

 \noi {\em Proof.} (i) Under given hypotheses, the uniform
 stability result in (ii) of Proposition \ref{Pr.91} implies  \cite[Ch.~1]{AMGV} that
 the $\o$-limit set of the asymptotically perturbed equation
 \ef{z3} is contained in that for the limit one \ef{z7}, which,
 under the given hypotheses,
 consists of the unique profile \ef{z12}.  \quad $\qed$





\subsection{BL-behaviour is ``almost" divergence free}

For the future convenience, we state again the asymptotic BL-behaviour: in the rescaled sense, by \ef{z10}, \ef{z12}, and
\ef{10},
 \be
 \label{9}
  \tex{
 \hat \vv(y,\t) = \cc_0(\t) \,  \ggg_0(y,\t)+... \, \whereA
  g_0^j(y,\t)= 1- \eee^{- \frac {\var(\t)}2\, {\rm dist}\{y, \pa \hat
  Q_0(\t)\}}.
  }
  \ee
 It is important that, by \ef{9}, the first term is ``a.e." an
 exponentially small perturbation of a divergence-free flow.
 Indeed, differentiating \ef{9} at any point staying away from
 the boundary by an arbitrarily small fixed ${\rm dist}\{\cdot\}=\d_0>0$, we have
 \be
 \label{9div}
  \tex{
  {\rm div}_y \, \hat \vv (y,\t) = O \big( c_0(\t)\var(\t) \eee^{- \frac
  {\var(\t)}2\,\d_0}\big) \to 0 \asA \t \to +\iy,
  }
  \ee
  provided that $c_0(\t)$ is not exponentially large (this does
  not happen, as we will show).
   In other words, not that surprisingly, the BL-expansion
 well
   keeps solenoidal features of originally divergence-free solutions $\hat \vv$ and, as
  customary, just makes an asymptotically (i.e., exponentially as in \ef{9div}) small
  perturbation of the ${\rm div}$. Thus, a somehow essential violation of the
  solenoidal property can happen only in an asymptotically small
  $O\big(\frac 1{\var(\t)}\big)$-neighbourhood of the boundary, which is negligible and  plays no role for $\t \gg 1$.

\section{Inner Region expansion: towards  ODEs regularity criterion}

 \label{S3.4}

\subsection{A standard semigroup approach leads to a more complicated problem}

Let us first perform necessary manipulations using a standard
semigroup approach. Applying to \ef{uu1} the projector ${\mathbb
P}$ yields
\be
  \label{uu1N}
   \tex{
  \hat \uu_t =
  \D \hat \uu - \frac {\pa \uu}{\pa \nn}\, \d_{S(t)}- {\mathbb P}\, (\hat \uu \cdot \n) \hat
  \uu.
  }
   \ee
Therefore, using the fundamental solution $b(x,t)$ of the heat
operator $D_t-\D$ with the rescaled kernel (the Gaussian) as in
\ef{BBB1} gives the following convolution representation of the
solution of \ef{uu1N}:
 \be
 \label{93}
  \tex{
  \hat \uu(t)= b(t)*\uu_0- \int\limits_0^t b(t-s)* \frac {\pa \uu}{\pa \nn}(s)\,
  \d_{S(s)}\, {\mathrm d}s - \int\limits_0^t b(t-s)*{\mathbb
  P}\, (\hat \uu(s) \cdot \n) \hat
  \uu(s)\, {\mathrm d}s.
 }
 \ee
 In particular, taking ${\rm div}$, we see that the second term on
 the right-hand side, containing  a surface integral, is div-free,
 since $\uu$ is.
Finally, sharply estimating the normal derivative therein and in
the third term from the core of BL-theory, the asymptotics \ef{9},
one can study the asymptotic behaviour of solutions, after using
the necessary scaling \ef{ll2}.

However, it turns out that this standard {\em integral} semigroup
approach leads to a more  complicated analysis, than a {\em
differential} one we will perform by using known spectral
properties of the rescaled operator $\BB^*$ involved in \ef{a3pp}.
Nevertheless, it is worth mentioning that such an approach can be
translated to the integral equation \ef{93}, with the clear
advantage of a more reliable rigorous justification by using
obviously smoother properties of solutions and, as a result, their
better uniform estimates in stronger metrics. It should be noted
that some principle difficulties cannot be avoided in such a way,
and, overall, technical technical questions  become  more dominant
for \ef{93}.

\subsection{Eigenfunction expansion: derivation of a 3D dynamical system}

Thus,  in the Inner Region, we deal with the original rescaled Cauchy
problem
 \ef{a3pp}.
Since, by construction, the extended solution orbit \ef{a1} is
uniformly bounded in $(L^2_{\rho^*}(\re^3))^3$, we can use the
converging in the mean (and uniformly on compact subsets in $y$)
eigenfunction expansion via the solenoidal Hermite polynomials
 as in \ef{bn12}:
 \be
 \label{a4}
  \tex{
  \hat \vv(y,\t)= \sum_{(\b)} c_\b(\t) \vv_\b^*(y).
   }
   \ee
   where we use the convention \ef{cb1} for the first vector $\cc_0(\t)$.
   Substituting \ef{a4} into \ef{a3} and using the orthonormality
   of these polynomials yield the following dynamical system for the
   expansion coefficients: for all $|\b| \ge 0$,
   \be
   \label{a5}
   \tex{
  \dot c_\b= \l_\b c_\b
-\langle
 \frac {\pa  \vv}{\pa \nn}\,
 \d_{\hat S(\t)}, \vv_\b \rangle
-  \eee^{-\frac \t 2} \,
   \langle {\mathbb P}(\hat \vv
 \cdot \n_y) \hat \vv H , \vv_\b \rangle ,
    }
    \ee
    where $\l_\b= -\frac {|\b|} 2$ are the real eigenvalues as in \ef{bbb1}.


 Recalling that eigenvalues in \ef{a5} satisfy
    $\l_\b \le - \frac 12$ for all $|\b| \ge 1$, it follows that
      we need to
     concentrate on the ``maximal" first Fourier generic pattern associated with
 the first constant Hermite polynomial $\vv_0^*$ in  \ef{Hp0},
 \be
 \label{a6}
 k=|\b|=0: \quad \l_0=0 \andA \vv_0^*(y) \equiv \eb=[1,1,1]^T.
  \ee
 The  normalized eigenfunction of the $L^2$-adjoint operator $\BB$
 is then
\be
 \label{a6F}
 \vv_0(y) = F(y)\, \eb,
  \ee
  where $F(y)$ is
 the Gaussian in
 (\ref{BBB2}).
 Actually, as follows from \ef{z10}, this corresponds
 to a naturally understood ``centre subspace behaviour"  for
 the operator $\BB^*$ in
  \ef{a3}:
  \be
  \label{a7}
  \hat \vv(y,\t) = \cc_0(\t)\eb + \ww^\bot(y,\t) \whereA \ww^\bot \in {\rm
  Span}\{\vv_\b^*, \,\, |\b| \ge 1\},
   \ee
   and $\ww^\bot(y,\t)$ is then negligible relative to $\cc_0(\t)$ for $\t \gg 1$.

   Fortunately, we actually do not need such a literal  using of
   those
   ``centre subspace issues" for  a difficult non-autonomous equation like \ef{a3pp},
    since the asymptotics of solutions
   \ef{a7} is directly  dictated by  Proposition \ref{Pr.Gen1}.
  In its turn,  this is an equivalent characterization of our class of generic
   patterns,
    and, in particular, the following holds:

\begin{proposition}
\label{Pr.Stable}
 Under the given above assumptions and hypotheses,  for the prescribed
class of generic
 solutions defined in Proposition $\ref{Pr.Gen1}$, \ef{a7} holds with
  $w^{\bot j}(\t)=o(|c_0^j(\t)|)$ as $\t \to +\iy$, and then the matching with the
 boundary layer behaviour in \ef{z2} requires
 \be
  \label{10}
   \tex{
   \frac {a_0^j(\t)}{\rho^j(s)} \to 1 \asA \t \to +\iy
   \LongA \rho^j(s) = a_0^j(\t)(1+o(1)), \quad j=1,2,3.
   }
   \ee
 \end{proposition}

 \noi {\em Proof.} This follows from the construction of the
 boundary layer by comparing the solution representations in \ef{z2}, \ef{z10}, and \ef{a7}.
 \quad $\qed$

 \ssk


\ssk

 Thus, the equation for $\cc_0(\t)$, with $\l_0=0$,  takes the form:
\be
   \label{a8}
   \tex{
 \dot \cc_0=-
\langle
 \frac {\pa  \vv}{\pa \nn}\,
 \d_{\hat S(\t)}, \vv_0 \rangle
- \eee^{-\frac \t 2}\,
   \langle {\mathbb P}(\hat \vv
 \cdot \n_y) \hat \vv H , \vv_0 \rangle,
    }
    \ee
where the first adjoint eigenfunction $\vv_0$ is as in \ef{a6F}.

We now need to return to BL-theory in Section \ref{S3}
establishing the boundary behaviour \ef{z2} for $\t \gg 1$, which
has the form \ef{9}.
 Then the convergence \ef{9}, which by a standard parabolic regularity is also true for the
 spatial derivatives, yields, in the natural rescaled sense,
  \be
  \label{11}
  \tex{
 \frac {\pa  \vv}{\pa \nn}
 =- \frac 12\, \cc_0(\t) \var(\t) \frac {y
 \cdot \nn(y)}{|y|}+...\, \whereA y \cdot \nn(y) <0 \,\,\, (\mbox{by convexity}).
 }
  \ee
  Therefore,
    the first (linear) term in \ef{a8}  asymptotically
reads
 \be
 \label{aa8}
  \tex{
-
\langle
 \frac {\pa  \vv}{\pa \nn}\,
 \d_{\hat S(\t)}, \vv_0 \rangle
  = \frac {\cc_0 \var(\t)} {2(4 \pi)^{3/2}}\,
   \int\limits_{\hat S(\t)} \frac {s
 \cdot \nn(s)}{|s|}\,   \eee^{-|s|^2/4} \, {\mathrm d}s + ...\,
 .
 }
 \ee


Due to the normalization in \ef{ph1}, we have that $a_i \le 1$, so
that the last term can be estimates above as: for any component
$j=1,2,3$,
 \be
 \label{aa81}
  \tex{
  \frac { \var(\t)} {2(4 \pi)^{3/2}}\,
   \int\limits_{\hat S(\t)} \frac {s
 \cdot \nn(s)}{|s|}\, \eee^{-|s|^2/4} \, {\mathrm d}s \le -
 \g_1
 \var^3(\t)\,  {\mathrm e}^{-\var^2(\t)/4}+...\, , \quad \g_1={\rm
 const.}>0.
 }
 \ee
 Recall that the extra factor $\var^2$ is obtained via integration
 over a closed {\em surface} in $\re^3$.

For the simple radial paraboloid shape in \ef{ph1N} and  \ef{QQ1},
with $a_i=1$, i.e., by \ef{rad11},
 \be
  \label{QQ1rad}
   \tex{
 \hat S(\t)=\pa \hat Q_0 \cap\{\t\}: \quad q_0(y) \equiv
 |y|= \var(\t), \andA
   \tilde S: \quad |z|=1,
   }
 \ee
\ef{aa81} presents a sharp asymptotics behaviour rather than an
estimate: for  $j=1,2,3$,
 \be
 \label{aa81rad}
  \tex{
  \frac {\cc_0 \var(\t)} {2(4 \pi)^{3/2}}\,
   \int\limits_{\hat S(\t)} \frac {s
 \cdot \nn(s)}{|s|}\, \eee^{-|s|^2/4} \, {\mathrm d}s = - \cc_0
 \g_1
 \var^3(\t)\,  {\mathrm e}^{-\var^2(\t)/4}+...\, .
 }
 \ee
  Naturally, it is possible to derive most sharp asymptotics for this radial case.

  \ssk

Finally, we need to estimate the last nonlinear quadratic term in
\ef{a8}: by \ef{9},
 \be
  \label{11N}
  \tex{
 \eee^{-\frac \t 2}\,\langle {\mathbb P}(\hat \vv
 \cdot \n_y) \hat \vv H , \vv_0 \rangle = \eee^{-\frac \t 2}\,
 (\cc_0 \cdot \eb)\cc_0  \, \frac 1{(4\pi)^{3/2}} \int\limits_{\re^3} {\mathbb P}[(
 \ggg_0
 \cdot \n_y) \ggg_0  H] \, \eee^{-|y|^2/4} \, {\mathrm d}y+ ... \, ,
  }
  \ee
 where the quadratic term $(\cc_0 \cdot \eb)\cc_0$ denotes the vector
 $(c_0^1+c_0^2+c_0^3)\cc_0 \in \re^3$.
 Fortunately, since this nonlinear term  enjoys having a fast
 decaying exponential factor $\eee^{-\frac \t 2}$, we do not need its better sharp
 estimates.
 We just need to show that, on the generic solutions obeying the
 BL-behaviour \ef{9}, this term always remains exponentially small,
 so does not play any role for the regularity conclusions.

We restrict to the radial case, though the asymptotic smallness
 of the convection can be similarly shown for more general convex
 paraboloids. Thus, using the representation of ${\mathbb P}(\hat \vv \cdot \n)\hat \vv $ given
 in \ef{HHH21} and using the change as in \ef{z1}, i.e., setting $y= \var(\t) z$, we
 have from \ef{9}, for $\t \gg 1$,
 \be
 \label{01}
  \begin{matrix}
 \eee^{- \frac \t 2}\,
 (\hat \vv \cdot \n)\hat \vv=  \eee^{- \frac \t 2}\,(\cc_0\cdot
 \eb)\cc_0\,
  \var\big[\big(1-\eee^{-\frac{\var^2}2\, d_z}\big)\eee^{-\frac{\var^2}2\,
  d_z}(\nn \cdot \eb)\quad
 \ssk\ssk\ssk\\
  - C_3 \var^2 \int \frac{z-\zeta}{|z-\zeta|^3} \, \eee^{-{\var^2}\,
  d_\zeta} \sum_{(i,j)}(\nn \cdot e_i)(\nn \cdot e_j) \, {\mathrm
  d} \zeta \big]+... \equiv
  J_1+J_2,
  \end{matrix}
  \ee
   where $d_z$ (and $d_\zeta$ in the integral) denotes the distance:
    $d_z= {\rm dist}\, \{z, \tilde S\}$.
 For such a rough estimate from above, one can omit the projector
 ${\mathbb P}$, using the fact that $\|{\mathbb P}\|=1$.

It is not difficult to see that, the $J_1$-term in \ef{01} is
leading to the following integral:
 \be
 \label{int11}
  \tex{
J_1 \sim \eee^{-\frac \t 2}\, (\cc_0 \cdot \eb)\cc_0 \var^4
\int\limits_{\{|z| \le 1\}}\eee^{-\frac {\var^2}4\,|z|^2}
\big(1-\eee^{-\frac{\var^2}2 \, d_z}\big) \eee^{-\frac{\var^2}2 \,
d_z} (\nn \cdot \eb)\, {\mathrm d} z, }
 \ee
 where, in the radial geometry, we may put $d_z=1-|z|$. Reducing
 the integral in \ef{int11} to a standard 1D one,
   it can be estimated above as follows: for some constant
   $\g_2>0$,
    \be
    \label{int21}
     \tex{
     |J_1(\t)| \le \g_2 \,  \eee^{-\frac \t 2}\,  | (\cc_0 \cdot \eb)\cc_0  | \, \var^4(\t)\, \eee^{-\frac{\var^2(\t)}4} \asA \t \to
     +\iy.
     }
     \ee

Consider the last  term  $J_2$ in \ef{01}. The corresponding upper
estimate is: for  $\g_{3,4}>0$,
 \be
 \label{int31}
 \begin{matrix}
 |J_2(\t)| \le \g_3 \eee^{-\frac \t 2}\,  |(\cc_0 \cdot \eb)\cc_0
 | \, \var^6\,
\big|\int\limits_{\{|z| \le 1\}}\eee^{-\frac {\var^2}4\,|z|^2}
\int
 \frac{z-\zeta}{|z-\zeta|^3} \,\eee ^{-\var^2 d_\zeta}
 \sum_{i,j}(\cdot)\, {\mathrm d} \zeta \, {\mathrm d} z\big|
  \qquad\quad
  \ssk\ssk\\
  \le \g_4 \eee^{-\frac \t 2}\, |  (\cc_0 \cdot \eb)\cc_0| \,
  \var^6(\t)\,
  \eee^{-\frac{\var^2(\t)}4} . \qquad\quad
  \end{matrix}
  \ee
  We do not guarantee that the $\var^6$ multiplier in the final estimate in \ef{int31}
  is any sharp  (as well as $\var^4$ in \ef{int21}), but it is
  sufficient for showing the convection neglect near the vertex.
  Indeed, any such very rough estimates (or omitting the projector ${\mathbb P}$ as we did above)
   cannot undermine the principal fact:  extra multipliers
   containing {\bf any} power of $\var(\t)$, being a {\em slow} growing function,  do
   not practically affect the exponentially decaying factor
   ${\mathrm e}^{-\t/2}$ as $\t \to +\iy$ in \ef{11N} and \ef{01}.
  Comparing with \ef{int11}, \ef{int21} yields that \ef{int31}
     supplies us with the leading term as $\t \to +\iy$.



\ssk

   Thus, bearing in mind all above assumptions and estimates for generic patterns,
   we obtain
 the following asymptotically approximate dynamical system for the first
expansion coefficients $\cc_0(\t)$: for  $\t \gg 1$,
 \be
 \label{12}
  \fbox{$
  \tex{
\dot \cc_0 \sim -  \g_1 \, \cc_0
 \var^3(\t) \, {\mathrm e}^{-\var^2(\t)/4}
+ \g_4 \,  \eee^{-\frac \t 2}\, (\cc_0 \cdot \eb)\cc_0   \,
\var^6(\t)\, \eee^{- \var^2(\t)/4}+... \, ,
 }
 $}
  \ee
  where we now omit all higher-order terms appeared via the above
hypotheses.
  The sign ``$\sim$" in \ef{12} means that, in the presentation of
  the influence of the nonlinear convection term, we  used the
  estimate \ef{int31}, rather than a sharp asymptotics. However,
  this estimate suffices   to declare that the convection
  term is negligible in the regularity analysis.

\subsection{3D regularity criterion}

Thus, according to \ef{Poi6}, the following conclusion holds:

\begin{theorem}
 \label{Th.Cr}
 Under the assumed above hypotheses and conditions, the vertex
 $(0,0)$ is regular for the NSEs problem \ef{1}, \ef{u01} in the class of generic solutions,
  iff
 \be
 \label{cr.11}
\fbox{$
  \mbox{the origin is globally asymptotically stable for
  the 3D dynamical system $(\ref{12})$},
   $}
   \ee
    i.e., any its solution satisfies:
  \be
  \label{CC1}
  \cc_0(\t) \to 0 \asA \t \to +\iy.
   \ee
   \end{theorem}


 \subsection{Two regularity conclusions}

 We begin with a simpler linear one.

 \ssk

\noi\underline{\sc 1. Linear Stokes problem}. As a first
consequence, we confirm that Petrovskii's criterion \ef{PCr1}
remains valid in the linear case. Recall that here, our analysis
do not include more difficult ``nonlinear" estimates used in
\ef{11N}.


\begin{theorem}
 \label{Th.Stokes1}
For the linear Stokes problem \ef{Stokes1}, under our hypotheses
on generic solutions, the regularity criterion of the vertex
$(0,0)$ is Petrovskii's one \ef{PCr1}.

  \end{theorem}

\noi{\em Proof.} Introducing the new time,
 \be
 \label{t1}
  \tex{
  \var^3(\t) \, \eee^{-\var^2(\t)/4} \, {\mathrm d}{\t}= {\mathrm d} s
  \LongA s= \int\limits_0^\t \var^3(\zeta) \, \eee^{-\var^2(\zeta)/4} \, {\mathrm
  d}{\zeta},
  }
  \ee
from \ef{12} (with no quadratic term), we obtain a linear diagonal
autonomous system
 \be
 \label{t2}
  \tex{
  \frac{{\mathrm d}}{{\mathrm d}s} \, \cc_0=-\g_1 \, \cc_0
   \LongA \cc_0(s)= \cc_0(0)\, \eee^{-\g_1 s} \to 0
  }
  \ee
 if and only if $s \to + \iy$ as $\t \to +\iy$, and the divergence
 in \ef{PCr1} follows. \quad $\qed$

 \ssk

 \ssk

\noi\underline{\sc 2. Navier--Stokes equations}. Performing the
change \ef{t1} for the full dynamical system \ef{12}  yields
 \be
 \label{12F}
  \tex{
\frac{{\mathrm d}}{{\mathrm d}s} \, \cc_0=-\g_1 \, \cc_0 + \g_4 \,
\eee^{-\frac \t 2}\, (\cc_0 \cdot \eb)\cc_0   \,   {\var^3(\t)}.
 }
  \ee
An elementary balancing of the linear and nonlinear term on the
right-hand side shows that the nonlinear one can be efficiently
involved into the {\em regular asymptotics} if it has at least an
exponential growth
 \be
 \label{t4}
  \tex{
 |\cc_0(\t)| \gg  \eee^{\frac \t 2}\,\frac 1{\var^3(\t)} \gg 1 \forA \t \gg 1,
 }
 \ee
 since, by assumptions, $\var(\t)$ is a slow growing function.
 Of course, \ef{t4} is impossible, since by the regularity assumption $\cc_0(\t)
 \to 0$.

 Note that, if $\vv(\t)$ gets vanishing as $\t \to +\iy$, then the
 same is true for the pressure via \ef{Poi4} (then \ef{Poi3} is a
 perturbation of the Laplacian), i.e., the pressure influence in
 the full equation (with no convection) is truly exponentially
 negligible.





\ssk

\ssk

 Thus,
 this is our final conclusion for the NSEs: under given hypothesis and for generic
 solutions\footnote{For the full Navier--Stokes models, in view of
 an essential difficulties to justify nonlinear convection
 estimates
 for calculating \ef{11N},
  we rather hesitate to state this even as a formal
 asymptotics; further mathematical research is necessary and is
 desirable.},
  {\em
the nonlinear convection term cannot affect the regularity of the
vertex $(0,0)$, so that Petrovskii's criterion \ef{PCr1} remains
true for the Navier--Stokes equations \ef{1}.}



\section{Final discussion: two key blow-up problems for the 3D NSEs}
 \label{SFi}

 It is now worth and  natural to look how the characteristic boundary
regularity analysis stands and fits in the wide area around
 the {\em Millennium Prize Problem
for the Clay Institute} (the MPPCI) on {\em global existence or
nonexistence of
 bounded smooth solutions for the NSEs \ef{1}}; see
Fefferman
 \cite{Feff00}. Both deal with  settings presented in a similar
 fashion:

 \ssk

 {\bf (${\mathbb B}{\mathbb R}$)} \underline{{\bf ${\mathbb B}$oundary point ${\mathbb R}$egularity}}: sharp ``blow-up" asymptotic expansions of
 possible solutions of \ef{1} {\em near a characteristic boundary point}  of $Q_0$. Then, mathematically speaking, after
 blow-up scaling \ef{ll2}, we arrive at an exponentially small
 perturbation of a solenoidal heat equation in \ef{S1}, which is
 convenient to write down in its full presentation
  \be
  \label{S1N}
  \fbox{$
   \tex{
    \vv _\t=
     \D \vv - \frac 12\, y \cdot \n \vv  - \fbox{$\eee^{-\frac \t 2}\,$} \,  {\mathbb P}(\vv
\cdot \n)\vv  \inB \hat Q_0.
  }
  $}
  \ee
 We observe here a crucial moment: in the vertex regularity study,
 the key nonlinear convection term gets an extra fast
 exponentially decaying factor $ \eee^{-\frac \t 2} $ (boxed in \ef{S1N}), that,
 indeed, essentially simplifies the matched asymptotic analysis.
 As we have seen, our goal was then simply to show that this
 nonlinear exponentially small perturbation does not affect the
 regularity of the vertex $(0,0)$ at all, so it becomes pure
 parabolic (Petrovskii's) one, i.e., governed as $\t \to +\iy$ by the rescaled
 solenoidal heat equation
  $$
  \tex{
\vv _\t=
     \D \vv - \frac 12\, y \cdot \n \vv   \inB \hat Q_0.
  }
  $$

 \ssk

 {\bf (${\mathbb I}{\mathbb R}$)} \underline{{\bf ${\mathbb I}$nterior point ${\mathbb R}$egularity (the MPPCI)}}: asymptotic expansion of possible
 solutions in {\em any interior point to check whether finite time
 blow-up in $L^\iy$ is possible or not}. Then,  a full
    Leray's  self-similar blow-up
   scaling must be performed\footnote{This Leray's blow-up self-similarity led him to a conjecture on {\em existence of a blow-up similarity solution}
as $t \to T^-$, and on {\em existence of  a self-similar extension} for $t>T$
(with $(T-t) \mapsto (t-T)$ in \ef{ll2Ler}); see some extra
details and a discussion in \cite[\S~2]{GalNSE}.}
   \cite{Ler34}:
   \be
   \label{ll2Ler}
   \tex{
   \uu(x,t)= \frac 1{\sqrt{T-t}}\, \vv (y,\t), \quad y = \frac
   x{\sqrt{T-t}}, \quad \t = - \ln(T-t),
}
 \ee
 where $T >0 $ is  the assumed finite blow-up time of the vector field $\uu(x,t)$ (used to be always  $T=0$
 beforehand).
 This
leads to  a different {\em autonomous} rescaled equation:
 \be
 \label{Full1}
  \fbox{$
  \tex{
  \vv _\t=
     \D \vv - \frac 12 \, y \cdot \n \vv - \frac 12 \, \vv  -
 \fbox{$
 \begin{matrix}
 \quad\\
 \quad
 \end{matrix}
  $}
    \,\,    {\mathbb P}(\vv
\cdot \n)\vv \inB \re^3 \times \re_+,
 }
 $}
 \ee
 where, unlike \ef{S1N}, {\em the box is empty}:  no exponentially decaying factor therein!



This is the main principal difference  between both key regularity
problems: in the rescaled equation \ef{Full1},
 the convection term {\bf is not accompanied} by an exponentially small
 time-factor as in \ef{S1N}. Therefore, it can play a vital role
 in creating a possible $L^\iy$-blow-up singularity. It turned out
 that a
 self-similar blow-up of Leray's type \ef{ll2Ler}, with a
 nontrivial
 profile
 $\vv=\vv(y)$ solving the stationary equation \ef{Full1},  does not exist:
  \be
  \label{F54}
   \tex{
\D \vv - \frac 12 \, y \cdot \n \vv - \frac 12 \, \vv  -   {\mathbb
P}(\vv \cdot \n)\vv=0 \inB \re^3, \,\,\, \vv \in L^2(\re^3) \LongA
\vv=0;
 }
 \ee
  see \cite{Nec96, Tsai98, Mil01, Hou07}. Therefore, it seems that a
 leading idea how to create a blow-up singularity for the NSEs is
 to deal with the so-called Type II blow-up solutions, i.e., those
 which violate a uniform similarity estimate for Type I solutions:
  \be
  \label{F55}
   \tex{
  \mbox{\bf Type I:} \quad \mbox{for some constant $C>0$, \, $|\uu(x,t)| \le  \frac C{\sqrt{T-t}}$ in $\re^3 \times
  (0,T)$},
  }
   \ee
   \be
   \label{F551}
    \tex{
    \mbox{i.e., {\bf Type II blow-up}:} \quad \limsup\limits_{t \to T^-} \sqrt{T-t}\, \sup\limits_{x \in
    \re^3}|\uu(x,t)|=+\iy.
 }
     \ee

In the rescaled variables \ef{ll2Ler}, Type II blow-up
   means looking for a global solution of \ef{Full1} that
  ``blows up" as $\t \to + \iy$. This problem is open, though some formal
  scenarios  of such a Type II blow-up in the NSEs have been discussed for a while;
   see, e.g.,  references and discussions
 in  \cite{GalNSE}.

 \ssk

 Overall, we  expect that the present study of simpler
 {\em boundary point regularity} governed by an exponentially
 perturbed rescaled equation \ef{S1N} is an inevitable and
 important step towards solving the main open {\em interior point regularity} blow-up problem for
 \ef{Full1}. Indeed, for \ef{S1N}, it turned out that the
 asymptotic behaviour near the vertex was spatially governed by the simplest
 {\em first solenoidal Hermite polynomial} (actually, a constant).
 While, in order to understand blow-up in  the ``MPPCI equation"
 \ef{Full1}, it seems that  an essential involvement of
 {\em all other} solenoidal Hermite polynomials $\{\vv_\b^*(y)\}$
 as eigenfunctions of the linear rescaled operator $\BB^*$
 (see Appendix A) should be detected first, {\em before} attacking
 a possible essentially nonlinear structure of a Type-II blow-up
 singularity for the NSEs.

\smallskip

\ssk

{\bf Acknowledgements.} The authors would like to thank
I.V.~Kamotski for useful discussions of various regularity issues
concerning  the NSEs and general PDE theory.



\begin{appendix}
\begin{small}
\section*{Appendix A: Hermitian  spectral theory of the linear rescaled
operator $\BB^*$:
 point spectrum and  solenoidal Hermite polynomials}
 \label{S2}
 \setcounter{section}{1}
\setcounter{equation}{0}


Thus, approaching the characteristic vertex $(0,0)$ in the blow-up
manner \ef{ll2}, one observes Hermite's operator $\BB^*$ as the
principal linear part of the rescaled equation \ef{S1}. Writing it
 the corresponding divergent form,
 \be
 \label{mm1}
 \tex{
    \BB^* \vv
 \equiv \frac 1 {\rho^*}\, \n \cdot (\rho^* \n \vv ),
 }
 \ee
 where the weight is
$\rho^*(y)={\mathrm e}^{-\frac{|y|^2}4}>0$,
 we observe that the actual rescaled evolution is now restricted to the
 weighted $L^2$-space $L^2_{\rho^*}(\re^3)$, with the
 exponentially decaying
  weight $\rho^*(y)$.
 Here, $ \BB^*$ is the  (``adjoint") Hermite operator with
 the point spectrum \cite[p.~48]{BS}
  \be
  \label{bbb1}
   \tex{
   \s(\BB^*)= \big\{ \l_k= - \frac{k}2, \quad
   k=|\b|=0,1,2,...\big\} \quad (\mbox{$\b$ is a multiindex in $\re^3$}),
    }
    \ee
    where each $\l_k$ has the multiplicity $\frac{(k+1)(k+2)}2$
     for $N=3$, or the binomial number $C_{N+k-1}^k$.
     The corresponding complete and closed set of eigenfunctions
     $\Phi^*=\{\psi_\b^*(y)\}$
is composed from separable Hermite polynomials. Note another
important  property of Hermite polynomials:
 \be
 \label{bn100}
 \forall \, \psi_\b^*, \quad \mbox{any derivative} \,\, D^\g  \psi_\b^* \quad \mbox{is also
 an eigenfunction with $k=|\b|-|\g| \ge 0$}.
  \ee
 Recall that \cite{BS}
  \be
  \label{bn10}
 \mbox{polynomial set} \,\,\, \Phi^* \,\,\, \mbox{is complete and closed in \,
  $L^2_{\rho^*}(\re^3)$}.
 \ee

Further spectral properties are convenient to demonstrate using
the linear operator $\BB$,
 \be
 \label{BBB1}
  \tex{
  \BB= \D + \frac 12\, y \cdot \n + \frac 32\, I \inB
  L^2_\rho(\re^3) \whereA \rho = \frac 1{\rho^*},
  }
  \ee
  which is adjoint to $\BB^*$ in the dual $L^2$-metric. It has the same point spectrum and the
  corresponding eigenfunctions are multiple of the same Hermite
  polynomials according to the well-known {\em generating
  formula}:
   \be
   \label{BBB2}
    \tex{
    \psi_\b(y)= \frac {(-1)^{|\b|}} {\sqrt{\b !}}\, D^\b F(y)
    \equiv \psi_\b^*(y) F(y)
    \whereA F(y)= \frac 1{(4\pi)^{3/2}} \, {\mathrm e}^{-|y|^2/4}
    }
    \ee
    is the rescaled kernel of the fundamental solutions of
    $D_t-\D$ in $\re^3 \times \re_+$.
Then, the  bi-orthonormality  holds:
  \be
  \label{bn1}
  \langle \psi_\b^*, \psi_\g \rangle=\d_{\b \g} \quad \mbox{for
  any} \quad \b, \, \g.
   \ee
Indeed, this dual metric can be also treated as that in
$L^2_{\rho^*}(\re^3)$ for the self-adjoint case, but we prefer to
keep `` $L^2$-dual" notations (for using also in non-symmetric
Burnett cases; cf. \ef{1B}).


\ssk

Obviously, one needs to consider eigenfunction expansions in the
solenoidal restriction
 \be
 \label{bn11}
\hat L^2_{\rho^*}(\re^3)= L^2_{\rho^*}(\re^3)^3\cap\{ \di \vv=0\}.
 \ee
Indeed, among the polynomials $\Phi^*=\{\psi_\b^*\}$ there are
many that well-suit the solenoidal fields. Namely, introducing the
eigenspaces
 $$
 \Phi_k^*= {\rm Span}\,\{\psi_\b^*, \,\, |\b|=k\}, \quad k \ge 1,
  $$
in view of \ef{bn100}, $\di$ plays a role of a ``shift operator"
in the sense that
 \be
 \label{sh1}
 \di: \Phi^{*3}_k \to \Phi^*_{k-1}.
  \ee

We next define  the corresponding solenoidal eigenspaces as
follows:
 \be
 \label{Sol1}
  \tex{
  {\mathcal S}_k^*= \{\vv^*=[v_1^*,v_2^*,v_3^*]^T: \quad {\rm div} \, \vv^*=0,
  \,\,v_i^* \in \Phi_k^*\}\whereA {\rm dim}\,\, {\mathcal S}_k^*=k(k+2);
  }
  \ee
see \cite{Gal02, Gal02A, Gallay06} and further references therein.

Actually, the paper \cite{Gal02} deals with {\em global}
asymptotics of NSEs solutions  as $t \to +\iy$, where the adjoint
operator $\BB$ in \ef{BBB1} occurs. Since $\BB$ is self-adjoint in
$L^2_\rho(\re^3)$, many results from \cite[Append.~A]{Gal02A} can
be applied to $\BB^*$.
 For a full collection,
 see \cite{Brand09, Brand07} for further asymptotic expansions and self-similar solutions.
In particular, this made it possible to construct therein fast
decaying solutions of the NSEs on each 1D stable manifolds with
the asymptotic behaviour\footnote{We present here only the first
term of expansion; as usual in dynamical system theory, other
terms in the case of ``resonance" can contain $\ln t$-factors
({\em q.v.} \cite{Ang88} for a typical PDE application); this
phenomenon was shown to exist
 for the NSEs in $\re^2$ \cite[p.~236]{Gal02A}.}
 \be
 \label{mn1}
  \tex{
 \uu_\b(x,t) \sim  \, t^{\l_k-\frac 12} \, \vv_\b \big( \frac x{\sqrt
 t}\big)+... \asA t \to \iy, \,\,\,\mbox{where} \,\,\, \vv_\b= {\vv_\b^*}F \in
 {\mathcal S}_k
   }
   \ee
 are solenoidal eigenfunctions of $\BB$.
  Namely, taking
    \be
    \label{kl1}
    \tex{
    \vv=[v_1,v_2,v_3]^T \in {\mathcal S}_k, \,\,\, v_i \in {\Phi}_k= {\rm Span}\,\big\{\psi_\b=
    \frac{(-1)^{|\b|}}{\sqrt{\b !}}\, D^\b F(y), \, |\b|=k \big\},
 }
  \ee
  where $F$ stands for the rescaled Gaussian in \ef{BBB2}, 
we have that
 \be
 \label{kl2}
  \tex{
 \di \vv = (v_1)_{y_1} + (v_2)_{y_2}+ (v_3)_{y_3}= \di (\vv^* F)
 \equiv (\di \vv^*) F - \frac 12\, y \cdot \vv^* \, F.
  }
  \ee
 This establishes  a one-to-one correspondence between
 solenoidal eigenfunction classes ${\mathcal S}_k^*$ in \ef{Sol1}
 for $\BB^*$ and ${\mathcal S}_k$ in \ef{mn1}
 for $\BB$; see \ef{Hp0}--\ef{HP2} below
 for the first  eigenfunctions
 $\vv_\b= \vv^*_\b F$. Therefore,
  ${\rm dim}\,\, {\mathcal
   S}_k=k(k+2)$, etc.;
see details and rather involved  proofs of the asymptotics
\ef{mn1} for $k=1$ and 2 in \cite{Gal02}.

\ssk

In particular,  those solenoidal Hermite polynomial eigenfunctions
of $\BB^*$ can be chosen as follows \cite[p.~2166-69]{Gal02A}
(the choice is obviously not unique, normalization constants are
omitted):
 \be
 \label{Hp0} \underline{\l_0=0:} \quad
 \vv_0^*=[1,1,1]^T=\eb \quad (\mbox{the first solenoidal Hermite polynomial}),
  \ee
 \be
 \label{HP1}
  \tex{
   \underline{\l_1=- \frac 12:} \quad
  \vv_{11}^*=\left[\begin{matrix}0\\-y_3\\y_2 \end{matrix} \right], \quad
\vv_{12}^*=\left[\begin{matrix}y_3\\0\\-y_1 \end{matrix} \right],
\quad \vv_{13}^*=\left[\begin{matrix}-y_2\\y_1\\0 \end{matrix}
\right]\,\,({\rm dim}\, {\mathcal S}_1^*=3);
 }
 \ee
 \be
 \label{HP2}
  \begin{matrix}
   \underline{\l_2=- 1:} \,\,\,
  \vv_{21}^*=\left[\begin{matrix}4-y_2^2-y_3^2\\y_1y_2\\-y_1y_3 \end{matrix} \right],
  \,\,
\vv_{22}^*=\left[\begin{matrix}y_1y_2\\4-y_1^2-y_3^2\\-y_2y_3
\end{matrix} \right], \,\,
\vv_{23}^*=\left[\begin{matrix}y_1y_3\\-y_2y_3\\4-y_1^2-y_2^2
\end{matrix} \right],
 \ssk\ssk\\
  \vv_{24}^*=-\left[\begin{matrix}0\\-y_1y_3\\y_1y_2 \end{matrix} \right], \quad
\vv_{25}^*=-\left[\begin{matrix}y_2y_3\\0\\-y_2y_1 \end{matrix}
\right], \qquad\qquad\qquad  \quad 
 \ssk\ssk\\
  \vv_{26}^*=\left[\begin{matrix}-y_2 y_3\\y_2y_3\\y_1^2-y_2^2 \end{matrix} \right],
  \,\,
\vv_{27}^*=\left[\begin{matrix}y_1y_2\\y_3^2-y_1^2\\-y_2y_3
\end{matrix} \right], \,\,
\vv_{28}^*=\left[\begin{matrix}y_2^2-y_3^2\\-y_1y_2\\y_1y_3
\end{matrix} \right]\,\,\,({\rm dim}\, {\mathcal S}_2^*=8),  \quad \mbox{etc.}
  \end{matrix}
 \ee

We need the following final conclusion.
By \ef{bn10}, the set of vectors $\Phi^{*3}$ is complete and
closed in\footnote{Note a standard result of functional analysis:
polynomials are  complete in any weighted $L^p$-space with an
exponentially decaying weight; see the analyticity argument in
Kolmogorov--Fomin \cite[p.~431]{KolF}.} $L^2_{\rho^*}(\re^3)^3$,
so that
 \be
 \label{bn12}
  \tex{
 \forall \, \vv \in L^2_{\rho^*}(\re^3)^3 \LongA \vv= \sum_{(\b)} c_\b
 \vv^*_\b, \quad \vv_\b^* \in \Phi^{*3}_k, \,\,\, k= |\b| \ge 0,
 }
  \ee
  where $xc_\b$ are scalars and, in a natural way,  the multiindex $\b$ arranges summation over  all  solenoidal
  Hermite polynomials.
   By obvious reasons, the only vector expansion  coefficient in \ef{bn12} is the first one, $\cc_0$, so, for convenience,  
  we will use the following vector notation:
  \be
  \label{cb1}
  \cc_0=[c_0^1,c_0^2,c_0^3]^T \LongA \cc_0
 \vv^*_0 \equiv
 [c_0^1 v_{0 1}^*,c_0^2v_{0 2}^*,c_0^3v_{0 3}^*]^T \quad (\vv_0^*=[1,1,1]^{T}).
 \ee

   It then follows
  from \ef{bn1}--\ef{sh1} that
   \be
  \label{bn10N}
 \mbox{polynomial set} \,\,\, \hat \Phi^*= \Phi^{*3}\cap\{\di \vv=0\} \,\,\,
 \mbox{is complete and closed in \,
  $\hat L^2_{\rho^*}(\re^3)$}.
 \ee
In what follows, we always assume that we  deal with ``solenoidal"
asymptotics involving eigenfunctions as in \ef{Sol1}.



\end{small}
\end{appendix}

\begin{appendix}
\begin{small}
\section*{Appendix B: Vertex regularity for Burnett equations}
 \label{SBur}
 \setcounter{section}{2}
\setcounter{equation}{0}

\subsection{Burnett equations in a hierarchy of hydrodynamic
models}

 The Burnett equations \ef{1B} appear as the {\em second
approximation} (the NSEs \ef{1} being the {\em first one}) of the
corresponding kinetic equations on the basis of Grad's method in
 Chapman--Enskog expansions for hydrodynamics.
Namely,   Grad's method applied to kinetic
 equations,  by expanding the kernel of the integral operators involved in terms
 of those with pointwise
 supports,
   yields, in addition to the classic operators of the Euler equations, other
  viscosity parts as follows:
 $$
  \mbox{$
  D_t \uu \equiv
  \uu_t +(\uu \cdot \n)\uu =  
  \sum\limits_{n=0}^\infty \e^{2n+1} \D^n(\mu_n \D \uu)+...= \e\big(\mu_0 \D \uu+ \e^2 \mu_1
  \D^2 \uu+...\big)+...\, ,
   $}
  $$
  where $\e>0$ is essentially the {\em Knudsen number} Kn;  see details in
  Rosenau's  regularization approach,
  \cite{Ros89}. In a full model, truncating such series at $n=0$ leads to the
  Navier--Stokes equations (\ref{1}) (with $\mu_0>0$), while $n=1$ is associated with the
    Burnett equations \ef{1B}.
Note  that Burnett-type equations, with a small parameter appeared
as higher-order viscosity approximations of the Navier--Stokes
equations, is an effective tool for proving existence of their
weak (``turbulent" in Leray's sense) solutions; see Lions'
monograph \cite[\S~6, Ch.~1]{JLi}.
Note that the
 ``Problem
 on blow-up/non-blow-up for
Burnett equations \ef{1B} at an interior point" starts from
dimensions $N = 7$: for $N \le 6$, there exists a unique global
smooth $L^2$-solution, \cite[\S~6]{G3KS}.

  The final finite-dimensional dynamical system
 for first Fourier coefficients of solutions to \ef{1B} is derived  similarly, but
 the regularity conclusions are shown to be more difficult and even rather obscure. The
 necessary spectral properties and vector solenoidal generalized Hermite
 polynomials as eigenfunctions of the rescaled non-self adjoint
 operator $\BB^*$ are introduced below in Appendix C.

More briefly, we now list below main steps of the regularity
analysis for \ef{1B}.

\subsection{First blow-up scaling: an exponentially perturbed parabolic equation}

  The first
   {\em blow-up} scaling in \ef{1B} is now
   \be
   \label{ll2B}
   \tex{
   \uu(x,t)=  \vv (y,\t), \quad y = \frac x{(-t)^{1/4}}, \quad \t=
   -\ln(-t): (-1,0) \to \re_+,
 }
 \ee
 which
 yields
 the following
rescaled equation:
 \be
  \label{S1B}
   \tex{
    \vv _\t= - \eee^{- \frac {3\t} 4} \n p + \BB^* \vv  -  \eee^{-\frac {3\t} 4}\,(\vv
\cdot \n)\vv  \inB \hat Q_0
 \whereA \BB^*=-\D^2   -
 \frac 14 \, y \cdot \n
  }
  \ee
 is the adjoint operator \ef{B1*} for $m=2$ with good
  spectral properties given in Appendix C.

 \subsection{Backward paraboloid}

Here,  $\pa Q_0$ is defined as follows:
 \be
 \label{ph1B}
  \tex{
 S(t)= \pa Q_0\cap\{t\}: \quad q_0(x) \equiv  \big( \sum_{i=1}^N a_i x_i^4\big)^{\frac 14}=
  (-t)^{\frac 14}\, \var(\t) \quad  (\sum a_i^4=1),
 }
  \ee
 with the same slow growing functions $\var(\t)$ as in \ef{ph1}.
The rescaled equation \ef{S1B} is then set in an expanding domain,
\be
 \label{ph1NB}
  \tex{
 \hat S(t)= \pa \hat Q_0\cap\{\t\}: \quad  q_0(y) \equiv \big(\sum_{i=1}^N a_i y_i^4\big)^{\frac 14}=
   \var(\t) \to +\iy
 \asA \t \to +\iy.
 }
  \ee

 \subsection{The Cauchy problem setting}

Extending to the Cauchy problem by using the variables \ef{a2}, we
use Green's second formula: for any $\chi \in C_0^\iy$,
 \be
  \label{Green2}
 \tex{
 \int_{\{q_0(y) \le \var(\t)\}} (\vv \D^2 \chi- \chi \D^2v) {\mathrm d}y= \int_{\hat S(t)}
\big(   \D \chi  \frac {\pa \vv }{\pa \nn} - \vv \frac {\pa \D
\chi}{\pa \nn}
 - \D \vv  \frac {\pa
\chi}{\pa \nn} + \chi  \frac {\pa \D \vv}{\pa \nn} \big) \,
{\mathrm d}s.
 }
 \ee
In view of the Dirichlet boundary conditions in \ef{1B}, in the
sense of distributions, \ef{Green2} reads:
 \be
 \label{i1BB}
  \tex{
 \D^2 \hat \vv= \D^2 \vv H + \big( \frac {\pa}{\pa \nn}\D \vv \big)\, \d_{\hat
 S(\t)} + \frac {\pa}{\pa \nn}(\D \vv \, \d_{\hat S(\t)}),
 }
 \ee
where densities of a single- and a double-layer
 potential now depend on the Laplacian $\D \vv$ instead of $\vv$ in \ef{a2}. Using the
 same harmonic pressure extension,
 we obtain
\be
 \label{a3ppB}
   \tex{
  \hat \vv_\t   =
   \BB^* \hat \vv - \eee^{-\frac {3\t} 4}\,{\mathbb P}\, (\hat \vv\cdot \n) \hat \vv
- {\mathbb P}\,(\frac {\pa}{\pa \nn}\D \vv) \, \d_{\hat
 S(\t)}
   -{\mathbb P}\, \frac {\pa}{\pa \nn}(\D \vv \, \d_{\hat S(\t)}).
}
  \ee
As for the NSEs with $m=1$, this
  problem is always locally well-posed, and is guaranteed to be
  globally well posed either for $N \le 6$ or for any sufficiently smooth initial data for $N \ge 7$.


 \subsection{Boundary layer and a perturbed rescaled equation}
   \label{S3.3B}

 Close to the lateral boundary of $Q_0$,
  the next rescaled variables are
  \be
  \label{z1B}
   \tex{
   z= \frac y{\var(\t)} \andA \hat \vv(y,\t)=\ww(z,\t).
   }
   \ee
This makes the corresponding rescaled paraboloid \ef{ph1NB} fixed:
 \be
 \label{QQ1B}
  \tex{
  \tilde S: \quad \big(\sum_{i=1}^N a_i z_i^4\big)^{\frac 14}=1.
   }
   \ee
The rescaled vector field $\ww$ satisfies a perturbed equation:
 \be
 \label{z1EqB}
  \begin{matrix}
 \ww_\t=-
 \frac 1{\var^4} \, \D^2_z \ww - \frac 14\, z \cdot \n_z \ww +
    \frac {\var'}\var \, z \cdot \n_z \ww - \frac 1 \var \,
      \eee^{-\frac {3\t} 4}\,{\mathbb P}(\ww \cdot \n_z)\ww
       \qquad
     \ssk\ssk\\
- \var^{N-4}\, {\mathbb P}\,\big(\frac {\pa}{\pa \nn}\D \vv\big)
\, \d_{\tilde
 S}
   - \var^{N-4}\,{\mathbb P}\, \frac {\pa}{\pa \nn}(\D \vv \, \d_{\tilde
   S}).
 \end{matrix}
    \ee
The BL-variables for a fixed point $z_0 \in \tilde S$ on the
boundary \ef{QQ1} are
  \be
  \label{z2B}
  \tex{
  \xi= \var^{\frac43 }(\t)(z_0-z), \quad \var^{\frac 43}(\t) {\mathrm d} \t={\mathrm d}s,
   \andA \ww(z,\t)= \rho(s) \ggg(\xi,s).
  }
  \ee
In this boundary layer,  we are looking for a generic
 pattern of the behaviour described by \ef{z1EqB} on compact subsets
 near the lateral boundary, satisfying
  \be
  \label{z44B}
  |\xi| = o\big(\var^{-\frac 43}(\t)\big) \to 0
   \LongA |z-z_0| = o\big(\var^{-\frac 83}(\t)\big) \to 0 \asA \t \to +
   \iy.
   \ee

Substituting \ef{z2B} into the PDE \ef{z1EqB} yields
 \be
 \label{z3B}
  \begin{matrix}
  \ggg_s= \AAA \ggg - \frac 14\,  \frac 1{\var^{4/3}} \, \xi \cdot \n_\xi \ggg - \frac
  {\var'_\t}{\var}\, \big(z_0- \frac \xi {\var^{4/3}}\big)\cdot \n_\xi \ggg \qquad\quad \ssk\ssk\\
- \frac 43\, \frac{\var'_\t}{ \var^{1/3}} \,
  \xi \cdot \n_\xi \ggg -  \frac {\rho'_s}{\rho} \,
  \ggg + \frac \rho{\var} \,  \eee^{-\frac {3\t} 4}\,{\mathbb P} (\ggg \cdot \n_\xi)\ggg
- \frac 1\rho \, \var^{N-4- \frac 43}\, {\mathbb P}\,\big(\frac
{\pa}{\pa \nn}\D \vv \big)\, \d_{\tilde
 S}
 \ssk\ssk\\
   - \frac 1 \rho\, \var^{N-4- \frac 43}\,{\mathbb P}\, \frac {\pa}{\pa \nn}(\D \vv \, \d_{\tilde
   S}),
  \quad \mbox{where} \quad
  \AAA \ggg=  -\D^2 \ggg + \frac 14\, z_0 \cdot \n_\xi \ggg.
    \end{matrix}
   \ee

 Again, in \ef{z3B}, we observe a perturbed linear uniformly parabolic
 equation. As usual, here one needs to
 check that all the perturbation terms are asymptotically small as $s \, (\mbox{or}\,\,\t) \to +\iy$,
 relative to the stationary autonomous operator $\AAA$.
 This is done similarly to  $m=1$ above.
  Overall,
the BL representation \ef{z2B} and \ef{g12}
 imply that \ef{z5} holds.
 Moreover,
asymptotically, the limit problem becomes one-dimensional,
depending on the space variable \ef{eta1}.


 We again pose the same asymptotic behaviour \ef{z5} at infinity.
 According to the scaling  \ef{z2B},
 let us fix a uniformly bounded rescaled orbit $\{\ggg(s),\,\, s>0\}$\footnote{As usual,
 the scaled function $\rho(s)$ remains unknown
 and to be determined by matching with the inner region behaviour.}. Then, by
   parabolic theory \cite{Fr, EidSys}, we can again pass to the
  limit in \ef{z3B}  along a subsequence $\{s_k\} \to +\iy$, removing small perturbations.
  Therefore,
   uniformly on compact subsets defined in
  \ef{z44B}, as $k \to \iy$,
   \be
   \label{z64B}
    \tex{
   \ggg(s_k+s) \to \hh(s) \whereA \hh_s=\AAA \hh, \quad \hh= \frac{\pa \hh}{\pa \nn}=0
   \,\,\,\mbox{at} \,\,\,\xi=0, \quad h^j|_{\xi=+\iy}=1.
    }
    \ee
The {limit equation} obtained from \ef{z3B},
 \be
  \label{z74B}
   \tex{
  \hh_s= \AAA \hh \equiv - \hh_{\eta\eta\eta\eta}+ \frac 14\, h_\eta
  }
   \ee
is again a standard linear parabolic PDE in $\ren \times \re_+$,
with a non self-adjoint operator
 $\AAA$, so  \ef{z74B} is not a gradient system in $L^2$.
 We then need to show that, in an appropriate weighted
 $L^2$-space if necessary and under the hypothesis \ef{z5},
 the stabilization holds, i.e.,
  the $\o$-limit set of the orbit $\{\hh(s)\}_{s>0}$ consists of a
  single
 equilibrium: as $s \to +\iy$,
  \be
  \label{z104}
   \left\{
   \begin{matrix}
  \hh(\xi,s) \to \ggg_0(\xi) \whereA \AAA \ggg_0=0 \,\,\, \mbox{for}
  \,\,\,\eta \in \re, \ssk\ssk \\
   \ggg_0=\ggg_0'=0 \quad \mbox{for} \quad \eta=0, \quad
  g^j_0(+\iy)=1.\qquad\,\,
  \end{matrix}
  \right.
   \ee
  This gives the unique solution of \ef{z104}
  (see \cite[\S~7]{GalPet2m} and \cite[\S~5]{GalMazf(u)}): for,
  e.g.,
  $z_0=[1,1,1]^T$,
 \be
 \label{z124}
  \tex{
  g^j_0(\xi)=1-{\mathrm e}^{- \frac \eta{2^{5/3}}} \,\big[ \cos \big( \frac
  {\sqrt{3} \, \eta}{2^{5/3}}\big) + \frac 1{\sqrt 3} \, \sin \big( \frac
  {\sqrt{3}\, \eta}{2^{5/3}}\big) \big], \quad j=1,2,...\, ,N.
  }
  \ee

It turns out that the limit problem \ef{z74B} possesses a number
of strong gradient and contractivity properties. Namely setting by
linearization
 \be
 \label{g65}
 \tex{
 \hh(s)=\ggg_0 + \ww(s) \LongA \ww_s= \AAA \ww \equiv - \ww_{\eta\eta\eta\eta} +
 \frac 14 \, \ww_\eta, \quad \ww=\ww_\eta=0 \,\,\, \mbox{at} \,\,\, \eta=0,
 }
 \ee
we arrive at the following (cf. Proposition \ref{Pr.91} for
$m=1$):

\begin{proposition}
 \label{Pr.Grad4}
 {\rm (i)} \ef{g65} is a gradient system in $L^2$, and

 {\rm (ii)} In the given class of solutions, the $\o$-limit set
 $\O_0$ of \ef{g65} consists of the origin only and  is
 uniformly stable.
 \end{proposition}

 \noi{\em Proof.} (i) One can see that \ef{g65} admits a
 monotone Lyapunov function obtained by multiplying by
 $\ww_{\eta\eta}$ in $L^2$:
  \be
  \label{g66}
 \tex{
   \frac 12\, \frac {\mathrm d}{{\mathrm d}s} \int ((w^j)_\eta)^2= -
   \int((w^j)_{\eta\eta\eta})^2 \le 0.
   }
    \ee
Hence, (ii) also follows. \quad $\qed$

\ssk

Thus, quite similar to the second-order case, under given
assumptions, we can pass to the limit $s \to +\iy$ along any
sequence in the perturbed gradient system \ef{z3B}. Then, again
similarly to $m=1$, the uniform stability of the stationary point
$g_0$ in the limit autonomous system \ef{z74B} in a suitable
metric
 guarantees that the asymptotically small perturbations do not
affect the omega-limit set; see \cite[Ch.~1]{AMGV}.
 However,
at this moment, we cannot avoid the following convention, which
for $m =2$ is much more key than for $m=1$.
 Actually, the convergence \ef{z64B} and
\ef{z104} for the perturbed dynamical system \ef{z3B} should be
considered as the main {\sc Hypothesis}, which characterizes the
class of generic patterns under consideration, and then the
normalization  \ef{z5} is its partial consequence.
 For
 bi-harmonic flows, a more clear characterization of this class of
 generic patterns is very difficult. It seems that a correct language
 of doing this (in fact, for both cases $m=1$ and $m \ge 2$)
  is to reinforce  the corresponding  centre subspace behaviour as in
 \ef{a7}.

Finally, we
 summarize these conclusions as follows:

\begin{proposition}
 \label{Pr.g04}


Under the given hypothesis and conditions, the
 problem $\ef{z3B}$ admits a family of solutions $($called  generic$)$
  satisfying $\ef{z104}$.

 \end{proposition}

Such a definition of generic patterns looks rather
non-constructive, which is unavoidable  for such higher-order
nonlocal PDEs.
However, \ef{z104} is expected to occur for ``almost all"
solutions.

\ssk

Thus, we stop further  discussions concerning the passage to the
limit $s \to +\iy$ in \ef{z3B}, which, as we have shown, under the
given hypotheses on the asymptotic smallness of perturbations
available, reduces to a linear stability analysis of the
nontrivial equilibrium $\ggg_0$ of the linear rescaled operator
$\AAA$ in \ef{z3B}. This has been resolved for a class of generic
solutions. More generally, we have to deal with solutions of
\ef{z3B} from the {\em stable subset} $\mathcal{W}_0$ of $\ggg_0$
within the prescribed perturbed equations \ef{z3B}. A clear,
constructive, and full identification of $\mathcal{W}_0$ is not
possible for such higher-order nonlocal perturbed parabolic
equations.

\ssk

 We summarize the conclusions as follows: {\em in what
follows, under the given hypothesis and conditions, we will deal
with a family of solutions $\mathcal{W}_0$ of \ef{z3B} $($called
generic$)$, for which \ef{z104} holds.}

\ssk





\subsection{Inner region analysis}
 \label{S3.44}

In the Inner Region,  the original rescaled problem
 \ef{a3ppB} occurs.
 For any  extended solution orbit \ef{a1}  uniformly bounded in
$L^2_{\rho^*}(\re)$, we  use  the eigenfunction expansion \ef{a4}
via the generalized solenoidal  Hermite polynomials \ef{t6}.
   Substituting \ef{a4} into \ef{a3ppB} and using the orthonormality
   property \ef{Ortog} yield a dynamical system: for any multiindex
   $\b$, with $|\b| \ge 0$, with the useful vector convention in $\ren$ of the same type as in \ef{cb1},
   \be
   \label{a54}
  \left\{
  \tex{
\dot c_\b= \l_\b c_\b - \langle {\mathbb P}\,\big(\frac
{\pa}{\pa \nn}\D \vv \big)\, \d_{\hat
 S(\t)}, \vv_\b \rangle   - \langle {\mathbb P}\, \frac {\pa}{\pa \nn}(\D \vv \, \d_{\hat S(\t)})
 , \vv_\b \rangle
 -  \eee^{-\frac {3\t} 4}\, \langle
{\mathbb P}(\hat \vv \cdot \n)\hat \vv, \vv_\b \rangle,
  }
    \right.
    \ee
    where $\l_\b= -\frac {|\b|} 4$ by \ef{spec1}, so that
    $\l_\b<0$ for any $|\b| \ge 1$.
As for $m=1$, one then needs to  concentrate on the first Fourier
generic patterns associated with
   the centre subspace for $\BB^*$ (cf. \ef{t62} for $N=3$)
 \be
 \label{a64}
 k=0: \quad \l_0=0 \andA \vv_0^*(y) =\eb= [1,1,...,1]^T, \quad
 \vv_0(y)=F(y)\eb.
  \ee
  This reflects
 another characterization of our class of generic patterns.
   The equation for the vector $\cc_0(\t)$ (see \ef{cb1}) then takes the form:
\be
   \label{a84}
 \tex{
\dot \cc_0=- \int\limits_{\hat S(\t)}\big(\frac {\pa}{\pa \nn}\D
\vv \big) \, \vv_0 \, {\mathrm d}s
 - \int\limits_{\hat S(\t)} \D \vv \,
  \frac {\pa}{\pa \nn} \vv_0 \, {\mathrm d}s
-  \eee^{-\frac {3\t} 4}\, \int\limits_{\ren} {\mathbb P}(\hat \vv
\cdot \n)\hat \vv  \, \vv_0 \, {\mathrm d}y.
     }
   \ee
 Note that first two terms on the right-hand side in \ef{a84} have
 a pure solenoidal parabolic (bi-harmonic) nature, while the only
 Navier--Stokes influence is presented by the last nonlinear term
 with an exponentially decaying factor.

Using next the boundary behaviour \ef{z104} with the 1D profile
\ef{z124} for $\t \gg 1$: in the rescaled sense, on the given
compact subsets, \ef{9} holds, with $\ggg_0$ given by \ef{z124},
where $\eta$ stands for the rescaled distance:
 \be
 \label{a84B}
  \tex{
   \eta=
\var^{\frac 13}(\t)\,{\rm dist} \{y, \hat S(\t)\}. }
 \ee
 Similar to \ef{9div}, such a BL-asymptotics is ``almost"
 solenoidal for $\t \gg 1$, i.e., perturbed by an exponentially small factor at
 any distance $\d_0>0$ from the boundary.
 By the matching of both Regions for such
 generic patterns, \ef{10} has to remain valid.

Performing, similar to \ef{11}--\ef{int31}, proper estimating  of
all the three terms on the right-hand side of \ef{a84} (the last,
the Navier--Stokes one can be again estimated rather roughly)
yields the following dynamical system for the first Fourier
coefficients:
 \be
 \label{124}
  \left\{
  \tex{
  \dot \cc_0= - \g_{11} \cc_0 \var^N \vv_0(\var)  -
  \g_{12} \cc_0 \var^{N- \frac 23} \vv'_0(\var)
   - \g_2 \,  \eee^{-\frac {3\t} 4}\, (\cc_0
  \cdot \eb)\cc_0 \, \var^{N+3}\, \vv_0(\var)+... \, ,
   }
    \right.
  \ee
 where, again, the first two terms are purely ``parabolic".
Similar to \ef{12}, we replace surface integrals by some
  ``average" values, actually assuming the radial dependence on
  $|y|$ with the rescaled surface $ \hat S(\t):
  \,\,\,\{|y|=\var(\t)\}$, as in \ef{QQ1rad}. As above, we do not
  guarantee that the multiplier $\var^{N+3}$ in the last term in
  \ef{124} is any optimal one (since the whole term will be shown to be
  negligible anyway).

Finally, using the expansion of the rescaled kernel given in
\ef{i5} and keeping the leading term only,
 the asymptotic dynamical
system reads
\be
 \label{124N}
  \left\{
  \begin{matrix}
  \dot \cc_0=  \frac {4 d_0}3 \g_{1} \cc_0 \var^{N-\d_0}(\t)\, \eee^{-d_0
  \var^{4/3}(\t)} [C_1 \sin(b_0 \var^{\frac 43}(\t))+C_2 \cos(b_0
  \var^{\frac 43}(\t))]\qquad\qquad\qquad \ssk\ssk\\
  \qquad - \g_2 \,  \eee^{-\frac {3\t} 4}\, (\cc_0
  \cdot \eb)\cc_0 \,  \var^{N+3-\d_0}(\t)\, \eee^{-d_0
  \var^{4/3}(\t)} [C_1 \sin(b_0 \var^{\frac 43}(\t))+C_2 \cos(b_0
  \var^{\frac 43}(\t))],
 \end{matrix}
  \right.
  \ee
where, as usual, $\g_{1,2} \in \ren$ are some constant vectors.

We then arrive at a typical and simple {\bf ODE regularity
criterion}: {\em under the given hypotheses and conventions of our
asymptotic analysis, the vertex $(0,0)$ is regular in the class of
generic solutions iff  any solution of the non-autonomous $3D$
dynamical system \ef{124N} vanishes as $\t \to +\iy$}, i.e., $0$
is globally asymptotically stable for \ef{124N}.

 \subsection{Two regularity conclusions}

{\bf 1.} The first regularity conclusion is straightforward: in
the absence of the convection term, i.e., for the {\em linear
fourth-order Stokes problem},
 \be
 \label{St4}
 \uu_t=-\n p-\D^2 \uu, \quad {\rm div}\, \uu=0 \inB Q_0,
  \ee
  with the Dirichlet boundary conditions as in \ef{1B},
 the vertex $(0,0)$ is regular provided that the {\em following integral
 diverges to $-\iy$} (cf. \ef{PCr1}):
  \be
  \label{pr1NN}
   \fbox{$
   \tex{
 \frac {4 d_0}3 \g_{1}\,
 \int\limits^{+\iy}
 \var^{N+\frac 13-\d_0}(s) \, \eee^{-d_0
  \var^{4/3}(s)} [C_1 \sin(b_0 \var^{\frac 43}(s))+C_2 \cos(b_0
  \var^{\frac 43}(s))] \, {\mathrm d}s=-\iy.
  }
  $}
  \ee
Since the function under the integrals is strongly oscillatory and, in general, is
of changing sign for $s \gg 1$, \ef{pr1NN} may require a special
procedure of ``oscillatory cut-off" of the given $\var(\t)$ to
delete a possible positive part of the divergent integral; see
\cite[\S~7]{GalPet2m}. In other words, to get the regularity
conclusion \ef{pr1NN}, the behaviour of $\var(\t)$ as $\t \to +\iy$
must be very carefully adjusted with the nonmonotone and
oscillatory behaviour of the rescaled kernel $F(\var(\t))$ of the
fundamental solution \ef{1.3R} of the bi-harmonic operator.

If \ef{pr1NN} fails, then the vertex is not regular. One can see
that the transition from a possible regularity (after an
oscillatory cut-off) to the guaranteed irregularity occurs at the
following ``critical" paraboloid with
 \be
 \label{var*}
 d_0 \var_*^{\frac 43}(\t)= \ln \t \LongA
 \var_*(\t)= d_0^{-\frac 34} (\ln \t)^{\frac 34} \equiv
 2^{\frac {11}4} \, 3^{- \frac 34} \, (\ln \t)^{\frac 34} \asA \t
 \to +\iy.
  \ee
  This is the fourth-order analogy of Petrovskii's function
  \ef{ph2}, so that the constant therein,
   \be
   \label{varC*}
    C_*= 2^{\frac {11}4} \, 3^{- \frac 34}
    \ee
    is optimal (similar to the ``2" in \ef{ph2}): replacing it by
    any larger one $C_*+\e$, with an $\e>0$, guarantees  convergence in \ef{pr1NN} and hence the
    irregularity of the  vertex $(0,0)$. Indeed, for such $\var(\t)$, the integral in
    \ef{pr1NN} simply {\em converges}, i.e.,  on almost all such centre subspace orbits,
 \be
 \label{ccc1}
  \cc_0(\t) \not \to 0 \asA \t \to +\iy.
   \ee

    \ssk

{\bf 2.}   Concerning the full nonlinear Burnett problem \ef{1B},
we
   arrive at the same conclusion as for the NSEs. Namely, it
   follows by balancing two terms on the right-hand side of
   \ef{124N}, that the last nonlinear one may be leading provided
   that
   \be
   \label{last1}
    \tex{
    |\cc_0(\t)| \gg  \eee^{\frac {3\t} 4}\, \frac 1{\var^{3}(\t)} \to +\iy \asA \t \to +\iy,
    }
     \ee
  which never happens in the case of the ``linear" regularity described by \ef{pr1NN}.
  Hence, the {\em convective term cannot change the vertex regularity},
  thus leading to the same regularity criterion.

\subsection{Nonexistence of similarity blow-up for Burnett equations: Type II
singularities are also needed}

For \ef{1B}, Leray's-type blow-up scaling \ef{ll2Ler} takes the
form
 \be
 \label{LL1}
 \tex{
  \uu(x,t)=(T-t)^{-\frac 34} \, \vv(y,\t), \quad y= \frac
  x{(T-t)^{1/4}}, \quad \t= - \ln(T-t),
  }
  \ee
  where $\vv(y,\t)$ solves the rescaled equation
   \be
   \label{LL2}
    \tex{
   \vv_\t= \BB^* \vv- \frac 34\, \vv- {\mathbb P}\,(\vv \cdot
   \n)\vv \inB \ren \times \re_+.
    }
    \ee
 Here $\BB^*$ is the linear rescaled operator \ef{B1*} with the
 known point spectrum and eigenfunctions being generalized Hermite
 polynomials.

 As for the NSEs \ef{1}
  (see Section \ref{SFi}),
  the first question  is whether
 a nontrivial Type I self-similar blow-up exists, i.e., whether
 a nontrivial stationary solution $\vv=\vv(y)$ of \ef{LL2} exists:
  \be
  \label{LL3}
   \tex{
-\D^2 \vv - \frac 14 \, y \cdot \n \vv - \frac 34\, \vv- {\mathbb
P}\,(\vv \cdot
   \n)\vv=0 \inB \ren, \quad \vv \in L^2(\ren).
   }
    \ee
 It is curious that a negative answer (i.e., similar for the NSEs in
  Section \ref{SFi}) can be obtained rather convincingly
  just by a local asymptotic analysis of the elliptic equation
 \ef{LL3}.
 As happens in practically all blow-up problems
for reaction-diffusion and other nonlinear PDEs (see  examples in,
e.g., \cite{Gal5types, AMGV, SGKM}),
  a ``generic" behaviour of its solutions
 as $z=|y| \to +\iy$ is governed by the leading lower-order linear terms,
 i.e., in the radial representation, this means that, for $z \gg
 1$,
  \be
  \label{LL4}
   \tex{
   - \frac 14 \, z\vv'_z - \frac 34 \, \vv +...=0
   \LongA \vv(z) \sim \frac {\bf C}{z^3} \asA z \to +\iy,
   }
   \ee
   where ${\bf C} \in \re^3$ is a constant vector. Of course, \ef{LL4}
   is just a rough radial estimate, so an extra ``angular separation" is
   necessary to produce all asymptotics like that at infinity.
However, \ef{LL4} is sufficient for a key negative conclusion: via
the local behaviour \ef{LL4}, for any ${\bf C} \ne 0$,
 \be
 \label{LL5}
 \tex{
  \frac {\bf C}{|y|^3}
  \in L^2(\{|y|>1\}) \quad \mbox{iff} \quad N<6.
  }
  \ee
  In other words, in the ``blow-up case"\footnote{For $N \le 6$,
  solutions of \ef{1B} do not blow-up in $L^\iy(\ren)$, \cite{G3KS}.} $N \ge 7$,
  blow-up cannot be of a self-similar (Type I) form \ef{LL1} with a nontrivial asymptotics \ef{LL4}.

Surely, this is not a proof of such a nonexistence, since the
 special single case ${\bf C}=0$ in \ef{LL4} has not been ruled
 out. Indeed, formally, it can happen that, for ${\bf C}=0$,
the similarity profile $\vv(y)$ solving \ef{LL3} may reach an
exponential decay at infinity (on derivation, see Appendix C)
 \be
 \label{LL51}
  \tex{
\vv(y) \sim {\bf C_1}\, \frac 1{ |y|} \, \eee^{-a_0 |y|^{4/3}} \whereA
a_0=3 \cdot 2^{-\frac 83} 
\andA {\bf C_1} \in \re^3.
 }
  \ee

\noi{\bf Example: a diversion to blow-up in a related semilinear
bi-harmonic flow.} As is known from similarity blow-up in
semilinear bi-harmonic equations such as
 \be
 \label{LL6}
 u_t=- \D^2 u + |u|^{p-1} u \inB \ren \times \re_+ \whereA p>1,
  \ee
such a behaviour is highly unlikely. To explain this, consider its
self-similar blowing-up solutions
 \be
 \label{LL7}
  \tex{
  u(x,t)=(T-t)^{-\frac 1{p-1}} \, v(y), \,\, y= \frac
  x{(T-t)^{1/4}} \,\, \Longrightarrow \,\, - \D^2 v - \frac 14\, y \cdot \n  v- \frac 1{p-1}
  \, v + |v|^{p-1} v=0.
  }
  \ee
 Checking the asymptotic behaviour as $y \to +\iy$, we again
 obtain the generic algebraic decay similar to \ef{LL4},  which is
 governed by two linear terms: as $z=|y| \to +\iy$,
  \be
  \label{LL8}
   \tex{
   - \frac 14\, z v'- \frac 1{p-1}
  \, v +...=0 \LongA
   v(y) \sim  C\, {|y|^{- \frac 4{p-1}}}.
    }
    \ee
    Similarly, for $C=0$, the behaviour gets exponentially decaying  (cf. \ef{LL51}):
    \be
    \label{LL9}
     \tex{
     v(y) \sim {C_1} |y|^\d \eee^{-a_0 |y|^{4/3}} \whereA \d=
     -\frac 23 \,\big(N- \frac 2{p-1}\big),
     }
     \ee
     with the same $a_0$ as in \ef{LL51}. See \cite[\S~2.3]{Gal5types} and Appendix C  for a derivation of  such
     two-scale WKBJ-asymptotics \ef{LL9} and \ef{LL51}.
Therefore, the asymptotic bundle of exponentially decaying
solutions \ef{LL9} of the ODE in \ef{LL7} contains a unique
parameter $C_1 \in \re$, i.e., it is {\bf one}-dimensional. Thus,
this is not enough to ``shoot" {\em two} symmetry conditions at
the origin:
 \be
 \label{LL10}
 v'(0)=v'''(0)=0,
  \ee
  so an extra parameter should be at hand, and this is $p$.
  As shown in \cite[\S~2]{Gal5types} by a careful
  numerical analysis of the ODE in \ef{LL7} for $N=1$, there exists a {\em unique
  value} of the exponent
   \be
   \label{LL11}
   p=p_\d=1.40... \, ,
   \ee
   for which \ef{LL7} admits a solution with the exponential decay
   \ef{LL9}, with some $C_1 \ne 0$.  Then, in the sense of bounded
   measures in $\re$, for $p=p_\d$,
   \be
   \label{LL12}
    \tex{
   |u(x,t)|^{\frac {p-1}4} \to D \d(x) \asA t \to T^-, \quad \mbox{with
   the
   constant} \,\,\, D= \int |v(y)|^{\frac{p-1}4}\, {\mathrm d}y>0.
   }
    \ee

\ssk

\noi{\bf Back to blow-up in Burnett equations.} We expect that a
similar phenomenon {\bf does not} exist for similarity blow-up in
the Burnett equations, i.e.,  exponentially decaying similarity
profiles \ef{LL51} do not exist. Recall that, unlike \ef{LL6}, the
equations \ef{1B} and \ef{LL3} {\bf do not} contain any free
parameter (like $p$ in \ef{LL6}), which could allow to get such a
solution at least for some its values. Of course, \ef{LL3} is a
system of three solenoidal fourth-order semilinear elliptic
equations, and a definite negative nonexistence conclusion is very
difficult to justify rigorously\footnote{An extra parameter may be
``hidden" in a kind of ``symmetry group" in the $\ren$-geometry
admitted by these PDEs (anyway, this looks  not that convincing).
 Overall, existence of a pure
self-similar blow-up for Burnett equations for $N=7$ is a too
simple way to settle this new ``fourth-order Millennium Problem",
and (at least one of the) authors would like to rule out such a
trivial solution of it.}.

\ssk

Overall, we arrive at the following plausible situation: {\em
similar to the NSEs \ef{1} in dimensions $N \ge 3$ $($see
discussion in Section \ref{SFi}$)$, blow-up in the Burnett
equations \ef{1B} in dimensions $N \ge 7$ cannot be self-similar
and requires constructing (or proving their nonexistence)
non-self-similar Type II blow-up singularities}\footnote{Here,
there occurs a {\em 4th-order Blow-up Problem for the Burnett
equations \ef{1B}} that may be much more difficult mathematically
than the  {\em Millennium Prize Problem for the Navier--Stokes
ones \ef{1}}. Of course, unlike the classic one in the actual
$\re^3$, a ``non-realistic" dimension $N=7$ makes the 4th-order
Problem less attractive for applications and for a general public,
but, mathematically, it can be even more fundamental for PDE
theory, since represents less understood features and principles
of interaction of a higher-order viscosity-diffusion operator with
a nonlinear convection one gathered in a nonlocal fashion.}.

\end{small}
\end{appendix}

\begin{appendix}
\begin{small}
\section*{Appendix C: Solenoidal Hermitian spectral theory for operator pair
$\{\BB, \,\BB^*\}$}
 \label{SAp2}
 \setcounter{section}{3}
\setcounter{equation}{0}

 We describe the necessary
 spectral properties of the linear $2m$th-order differential
 operator in $\ren$ ($m=2$ for the Burnett  equations \ef{1B})
\begin{equation}
 \label{B1*}
  \tex{
 \BB^* = (-1)^{m+1} \D^{m}_y - \frac {1}{2m}
 \, y \cdot \n_y,
 }
 \end{equation}
and of its $L^2$-adjoint $\BB$
 given
by
  \begin{equation}
 \label{B1}
  \tex{
 \BB = (-1)^{m+1} \D_y^{m} + \frac {1}{2m}
 \, y \cdot \n_y +  \frac{N}{2m}\, I.
 }
 \end{equation}
 As we have seen, for $m=1$, \ef{B1*} and \ef{B1}  are classic Hermite self-adjoint  operators
 with completely known spectral properties, \cite[p.~48]{BS}.
 For any $m \ge 2$, both operators \ef{B1*} and \ef{B1},
 though looking very similar to those for $m=1$,
 {\em are not symmetric}  and do not admit a self-adjoint
extension, so we follow \cite{Eg4} in presenting  spectral theory.


\subsection{Fundamental solution, rescaled kernel, and first estimates}
 \label{Sect3}

  The fundamental solution $b(x,t)$  of the linear poly-harmonic parabolic equation
  \be
  \label{Lineq}
  u_t = - (-\D)^m u \inA
  \ee
 takes the standard similarity form
  \be
  \label{1.3R}
   \tex{
   b(x,t) = t^{-\frac N{2m}}F(y), \quad y= \frac  x{t^{1/{2m}}}.
    }
 \ee
 The rescaled kernel $F$ is the unique radial solution of the elliptic
 equation
  \begin{equation}
\label{ODEf}
 {\bf B} F \equiv -(-\Delta )^m F + \textstyle{\frac 1{2m}}\, y \cdot
\nabla F + \textstyle{\frac N{2m}} \,F = 0
 \,\,\,\, {\rm in} \,\, \ren,  \quad \mbox{with} \,\,\, \textstyle{\int F =
 1.}
\end{equation}
 For  $m \ge 2$, the rescaled kernel function $F(|y|)$ is   oscillatory as $|y| \to \infty$ and
satisfies
 \cite{EidSys, Fedor}
\begin{equation}
\label{es11} 
 |F(y)| < D\,\,  {\mathrm e}^{-d_0|y|^{\alpha}}
\,\,\,{\rm in} \,\,\, \ren, \quad \mbox{where} \,\,\,
\a=\textstyle{ \frac {2m}{2m-1}} \in (1,2),
\end{equation}
for some positive constants $D$ and $d_0$ depending on $m$ and
$N$.

\subsection{Some constants}

 As we have seen,   the rescaled kernel $F(y)$
satisfies \ef{es11}, where $d_0$ admits an explicit expression;
see below.
 Such optimal exponential estimates of the fundamental solutions
 of higher-order parabolic equations are well-known and were first
 obtained by Evgrafov--Postnikov (1970) and Tintarev (1982); see
 Barbatis \cite{Barb, Barb04} for key references.

As a crucial  issue for the  boundary point regularity study, we
will need a sharper, than given by \ef{es11}, asymptotic behaviour
of the rescaled kernel $F(y)$ as $y \to +\iy$. To get that, we
keep four  leading terms in
\ef{ODEf} and
obtain, in terms of the radial variable $y \mapsto |y|>0$:
 \be
 \label{i1}
 \tex{
 (-1)^{m+1} \big[F^{(2m)} + m \frac{N-1}y \, F^{(2m-1)}+...\big]
  + \frac 1{2m} \, y F' + \frac N{2m} \, F=0 \forA y \gg 1.
 }
 \ee
 Using standard classic WKBJ asymptotics, we substitute into \ef{i1}
 the function
  \be
  \label{i2}
  F(y) = y^{-\d_0} \, {\mathrm e}^{a y^\a}+... \asA y \to + \iy,
   \ee
   exhibiting two scales.
   Balancing two leading terms
 gives the algebraic equation for $a$ and  $\d_0$:
 \be
 \label{i3}
  \tex{
 (-1)^m (\a a)^{2m-1}= \frac 1{2m} \andA \d_0=  \frac{m(2N-1)-N}{2m-1}>0\,
 .
 }
  \ee



By construction, one needs to get the root $a$ of \ef{i3} with the
maximal ${\rm Re}\, a<0$. This yields (see e.g., \cite{Barb,
Barb04} and \cite[p.~141]{GSVR})
 \be
 \label{i4}
  \tex{
 a= \frac{2m-1}{(2m)^\a} \big[-\sin\big( \frac{\pi}{2(2m-1)}\big) +
 \ii \cos\big( \frac{\pi}{2(2m-1)}\big)\big] \equiv -d_0 + \ii b_0
 \quad (d_0>0).
 }
 \ee
Finally, this gives the following double-scale asymptotic of the
kernel:
 \be
 \label{i5}
  \tex{
  F(y) =
   y^{-\d_0} \, {\mathrm e}^{-d_0 y^\a} \big[ C_1 \sin (b_0 y^\a)+
   C_2 \cos (b_0 y^\a)\big]+... \asA y=|y| \to + \iy ,
   }
   \ee
 where $C_{1,2}$ are real constants, $|C_1|+|C_2| \not = 0$.
 In \ef{i5}, we present the first two leading terms from the
 $m$-dimensional bundle of exponentially decaying asymptotics.

In particular, for the Burnett equations \ef{1B} in $\re^3$, we
have
 \be
 \label{i6}
  \tex{
  m=2, \,\,N=3: \quad \a= \frac 43,  \quad d_0=3 \cdot 2^{-\frac{11}3},
 \quad b_0=3^{\frac 32} \cdot 2^{-\frac{11}3},
  \andA \d_0= \frac
  73.
   }
   \ee

\subsection{The discrete real spectrum and eigenfunctions of
 $\BB$}



 For $m \ge 2$,
 $\BB$ is considered  in the weighted space $L^2_\rho(\ren)$ with the
exponentially growing weight function
 \be
  \label{rho44}
  \rho(y) = {\mathrm e}^{a |y|^\a}>0 \quad {\rm in} \,\,\, \ren,
 \ee
  where $a \in (0,  2d_0)$ is a fixed
constant.
 We next
introduce a standard  Hilbert (a weighted Sobolev) space of
functions $H^{2m}_{\rho}(\ren)$ with the inner product and the
induced  norm
\[
 \tex{
 \langle v,w \rangle_{\rho} = \int\limits_{\ren} \rho(y) \sum\limits_{k=0}^{2m}
 D^{k}_y
 v(y) \, \overline {D^{k}_y w(y)} \,{\mathrm d} y \andA
\|v\|^2_{\rho} = \int\limits_{\ren} \rho(y) \sum\limits_{k=0}^{2m}
|D^{k}_y
 v(y)|^2 \, {\mathrm d} y.
  }
\]
Then $H^{2m}_{\rho}(\ren) \subset L^2_{\rho}(\ren) \subset
L^2(\ren)$, and  $\BB$ is a bounded linear operator from $
H^{2m}_{\rho}(\ren)$ to $ L^2_{\rho}(\ren)$. Key spectral
properties of the operator $\BB$ are as follows \cite{Eg4}:

\begin{lemma}
\label{lemspec}
 {\rm (i)}  The spectrum of $\BB$
comprises real simple eigenvalues only,
 \begin{equation}
\label{spec1}
 \tex{
 \sigma(\BB)=
\big\{\lambda_\b = -\frac k{2m}, \,\, k= |\b|= 0,1,2,...\big\}.
 }
\end{equation}
 {\rm (ii)} The eigenfunctions $\psi_\b(y)$ are given by
\begin{equation}
\label{eigen} \psi_\beta(y) =\textstyle{\frac{(-1)^{|\b|}}{\sqrt
{\b !}}} \, D^\beta F(y), \quad \mbox{for any} \,\,\,|\b|=k \ge 0.
\end{equation}

\noi{\rm (iii)} Eigenfunction subset \ef{spec1} is complete
 in $L^2({\re})$ and in $L^2_{\rho}({\re})$.

\noi {\rm (iv)} The resolvent $(\BB-\lambda I)^{-1}$
for  $\lambda \not \in \sigma(\BB)$ is a compact integral operator
in $L^2_{\rho}(\ren)$.
\end{lemma}


By Lemma \ref{lemspec}, the   centre and stable subspaces of $\BB$
are given by
 \be
 \label{centr1}
E^c = {\rm Span}\{\psi_0= F\} \andA E^s = {\rm Span}\{\psi_\b, \,
|\b|>0\}.
 \ee

\subsection{Polynomial eigenfunctions of the operator $\BB^*$}

Consider the operator (\ref{B1*}) in the weighted space
$L^2_{\rho^*}(\ren)$, where $\langle \cdot, \cdot
\rangle_{\rho^*}$ and $\|\cdot\|_{\rho^*}$ being the inner product
and the norm,
 with the   ``adjoint"  exponentially decaying weight
function
  \begin{equation}
\label{rho2}
 \tex{
 \rho^*(y) \equiv \frac  1 {\rho(y)} = {\mathrm e}^{-a|y|^{\a}} > 0.
 }
\end{equation}
 We ascribe to $\BB^*$ the domain
 $H^{2m}_{\rho^*}(\ren)$, which is dense in $L^2_{\rho^*}(\ren)$, and
 then
$$
 \BB^*: \,\, H^{2m}_{\rho^*}(\ren) \to L^2_{\rho^*}(\ren)
$$
 is a bounded linear operator.  $\BB$ is adjoint
 to $\BB^*$ in the usual sense: denoting by $\langle \cdot,\cdot \rangle $  the
inner product in the dual space $L^2(\ren)$, we have
  \begin{equation}
 \label{Badj1}
 \langle \BB v, w \rangle =  \langle v, \BB^* w \rangle
 \quad \mbox{for any} \,\,\, v \in H^{2m}_\rho(\ren) \andA
 w \in H^{2m}_{\rho^*}(\ren).
 \end{equation}
The eigenfunctions of $\BB^*$ take a particularly simple finite
polynomial form and are as follows:


 \begin{lemma}
\label{lemSpec2}
 {\rm (i)} $ \sigma(\BB^*)=\s(\BB)$.

 \noi{\rm (ii)} The eigenfunctions  $\psi^*_\b(y)$ of $\BB^*$ are generalized Hermite
 polynomials of degree $|\b|$ given by
 \begin{equation}
 \label{psi**1}
  \tex{
 \psi_\b^*(y) = \frac 1{\sqrt{\beta !}}
 \Big[ y^\b + \sum_{j=1}^{[|\b|/2m]} \frac 1{j !}(-\Delta)^{m j} y^\b
 \Big] \quad \mbox{for any} \quad \b.
 }
 \end{equation}

 \noi{\rm (iii)} Eigenfunction subset \ef{psi**1} is complete   in $L^2_{\rho^*}(\ren)$.

  \noi {\rm (iv)}
$\BB^*$ has a compact resolvent $(\BB^*-\lambda I)^{-1}$ in
$L^2_{\rho^*}(\ren)$ for $\lambda \not \in \sigma(\BB^*)$.

\noi{\rm (v)}  The bi-orthonormality of the bases $\{\psi_\b\}$
and $\{\psi_\g^*\}$ holds in the dual $L^2$-metric:
 \be
 \label{Ortog}
\langle \psi_\b, \psi_\g^* \rangle = \d_{\b\g} \quad \mbox{for
any} \quad \b,\, \g.
 \ee
\end{lemma}





\ssk

\noi{\bf Remark on closure.} This is an important issue for using
eigenfunction expansions of solutions.
 Firstly,
  in the self-adjoint case $m=1$, the sets of eigenfunctions are closed in the
corresponding spaces, \cite{BS} (and we have used this in our
previous NSEs study).

Secondly, for $m \ge 2$, one needs some extra details. Namely,
using (\ref{Ortog}), we can introduce the subspaces of
eigenfunction expansions and begin with the operator $\BB$. We
denote by $\tilde L^2_\rho$ the subspace of eigenfunction
expansions $v= \sum c_\b \psi_\b$ with coefficients $c_\b =
\langle v, \psi^* \rangle$ defined as the closure of the finite
sums $\{\sum_{|\b| \le M} c_\b \psi_\b\}$ in the norm of
$L^2_\rho$. Similarly, for the adjoint operator $\BB^*$, we define
the subspace $\tilde L^2_{\rho^*} \subseteq L^2_{\rho^*}$. Note
that since the operators are not self-adjoint and the
eigenfunction subsets are not orthonormal, in general, these
subspaces can be different from
$
 L^2_{\rho}$ and $L^2_{\rho^*}$, and particularly
the equality is guaranteed in the self-adjoint case $m=1$,
$a=\frac 1 4$.

Thus, for $m \ge 2$, in the above subspaces obtained via a
suitable closure, we can apply standard eigenfunction expansion
techniques as in the classic self-adjoint case $m=1$.


%

\subsection{Solenoidal Hermite polynomials}

The vector solenoidal Hermite polynomials are constructed from
\ef{psi**1} in a manner similar to that for $m=1$; cf
\ef{Hp0}--\ef{HP2}. Namely, given a vector polynomial
 \be
 \label{t6}
 \vv_\b^*=[\psi_{\b_1}^*, \psi_{\b_2}^*,..., \psi_{\b_N}^*]^T
 \whereA |\b_1|=|\b_2|=...=|\b_N|=|\b|,
  \ee
 it gets solenoidal provided that
  \be
  \label{t61}
   \tex{
   {\rm div}\, \vv_\b^* \equiv \sum\limits_{i=1}^N (\psi_{\b_i})_{y_i}
   =0.
   }
   \ee

For instance, for the Burnett case $m=2$ and $N=3$, some pairs are (not all
linearly independent eigenfunctions  are presented, normalization constants
are omitted):
 \be
 \label{t62}
 \l_0=0: \quad \vv_0^*=[1,1,1]^T,
  \ee
  \be
  \label{t63}
   \tex{
   \l_1=-\frac 14: \vv_{11}^*=[y_2, -y_3, y_2]^T, \quad
   \vv_{12}^*=[y_3,y_3,-y_1]^T, \quad \vv_{13}^*=[-y_2,y_1,y_1]^T,
   }
    \ee
 \be
 \label{t64}
  \tex{
 \l_2=- \frac 12: \quad \vv_{21}^*=[-y_1^2-y_3^2,y_1y_2,y_1y_3]^T, \,\,
 \vv_{22}^*=[y_1y_2,-y_2^2-y_3^2,y_2y_3]^T, \,\,\, \mbox{etc.}
 }
  \ee
  \be \label{t65}
   \tex{
  \l_3=- \frac 34: \quad \vv_{31}^*=[y_2^3,y_3^3,y_1^3], \quad
   \vv_{32}^*=
  [y_1y_2^2, y_2y_1^2,-y_3(y_1^2+y_2^2)], \,\,\,\mbox{etc.}
  }
   \ee
    \be
    \label{t66}
    \l_4=-1: \quad \vv_{41}^*=[y_2^4+4!,y_3^4+4!,y_1^4+4!]^T, \,\,\,
     \vv_{42}^*=[y_1y_2^3,y_2y_1^3,-y_3(y_1^3+y_2^3)]^T, \quad
     \mbox{etc.}
     \ee
As in the self-adjoint case $m=1$, some technical efforts are
necessary toward  completeness/closure of  generalized solenoidal
Hermite polynomials in suitable spaces. We omit details.
\end{small}
\end{appendix}

\end{document}